%
%
%
%
%
\def\fileversion{1.14c}
\def\filedate{2006/04/29}
\csname PSTricksLoaded\endcsname
\PstAtCode{\the\catcode`\@}
\catcode`\@=11\relax
\expandafter\ifx\csname @latexerr\endcsname\relax
\long\def\@ifundefined#1#2#3{\expandafter\ifx\csname
  #1\endcsname\relax#2\else#3\fi}
\def\@namedef#1{\expandafter\def\csname #1\endcsname}
\def\@nameuse#1{\csname #1\endcsname}
\def\@eha{%
  Your command was ignored.^^J
  Type \space I <command> <return> \space to replace
  it with another command,^^J
  or \space <return> \space to continue without it.}
\def\@spaces{\space\space\space\space}
\def\typeout#1{\immediate\write\@unused{#1}}
\alloc@7\write\chardef\sixt@@n\@unused
\def\@empty{}
\def\@gobble#1{}
\def\@nnil{\@nil}
\def\@ifnextchar#1#2#3{%
\let\@tempe#1\def\@tempa{#2}\def\@tempb{#3}\futurelet\@tempc\@ifnch}
\def\@ifnch{%
  \ifx\@tempc\@sptoken
    \let\@tempd\@xifnch
  \else
    \ifx\@tempc\@tempe \let\@tempd\@tempa \else \let\@tempd\@tempb \fi
  \fi
  \@tempd}
\begingroup
\def\:{\global\let\@sptoken= } \:
\def\:{\@xifnch} \expandafter\gdef\: {\futurelet\@tempc\@ifnch}
\endgroup
\fi
\typeout{`PSTricks' v\fileversion\space\space <\filedate> (tvz)}
\def\@pstrickserr#1#2{%
  \begingroup
  \newlinechar`\^^J
  \edef\pst@tempc{#2}%
  \expandafter\errhelp\expandafter{\pst@tempc}%
  \typeout{%
    PSTricks error. \space See User's Guide for further information.^^J
    \@spaces\@spaces\@spaces\@spaces
    Type \space H <return> \space for immediate help.}%
  \errmessage{#1}%
  \endgroup}
\def\@ehpa{%
  Your command was ignored. Default value substituted.^^J
  Type \space <return> \space to procede.}
\def\@ehpb{%
  Your command was ignored. Will recover best I can.^^J
  Type \space <return> \space to procede.}
\def\@ehpc{%
  You better fix this before proceding.^^J
  See the PSTricks User's Guide or ask your system administrator for help.^^J
  Type \space X <return> \space to quit.}
\def\pst@misplaced#1{\@pstrickserr{Misplaced \string#1 command}\@ehpb}
\newdimen\pst@dima
\newdimen\pst@dimb
\newdimen\pst@dimc
\newdimen\pst@dimd
\newdimen\pst@dimg
\newdimen\pst@dimh
\newbox\pst@hbox
\newbox\pst@boxg
\newcount\pst@cnta
\newcount\pst@cntb
\newcount\pst@cntc
\newcount\pst@cntd
\newcount\pst@cntg
\newcount\pst@cnth
\newif\if@pst
\newtoks\pst@toks
\newif\if@star
\def\pst@ifstar#1{%
  \@ifnextchar*{\@startrue\def\next*{#1}\next}{\@starfalse#1}}
\def\pst@expandafter#1#2{%
  \def\next{#1}%
  \edef\@tempa{#2}%
  \ifx\@tempa\@empty
    \@pstrickserr{Unexpected empty argument!}\@ehpb
    \def\@tempa{\@empty}%
  \fi
  \expandafter\next\@tempa}
\def\pst@dimtonum#1#2{\edef#2{\pst@@dimtonum#1}}
\def\pst@@dimtonum#1{\expandafter\pst@@@dimtonum\the#1}
{\catcode`\p=12 \catcode`\t=12 \global\@namedef{pst@@@dimtonum}#1pt{#1}}
\def\pst@pyth#1#2#3{
  \begingroup
    \pst@dima=#1\relax
    \ifnum\pst@dima<\z@\pst@dima=-\pst@dima\fi  
    \pst@dimb=#2\relax
    \ifnum\pst@dimb<\z@\pst@dimb=-\pst@dimb\fi  
    \advance\pst@dimb\pst@dima         
    \ifnum\pst@dimb=\z@
      \global\pst@dimg=\z@             
    \else
      \multiply\pst@dima 8\relax              
      \pst@@divide\pst@dima\pst@dimb     
      \advance\pst@dimg -4pt            
      \multiply\pst@dimg 2
      \pst@dimtonum\pst@dimg\pst@tempa
      \pst@dima=\pst@tempa\pst@dimg           
      \advance\pst@dima 64pt         
      \divide\pst@dima 2\relax                      
      \pst@dimd=7pt                
      \pst@@pyth\pst@@pyth\pst@@pyth            
      \pst@dimtonum\pst@dimd\pst@tempa
      \pst@dimg=\pst@tempa\pst@dimb
      \global\divide\pst@dimg 8             
    \fi
  \endgroup
  #3=\pst@dimg}
\def\pst@@pyth{
  \pst@@divide\pst@dima\pst@dimd
  \advance\pst@dimd\pst@dimg
  \divide\pst@dimd 2\relax}
%
%
\def\pst@divide#1#2#3{%
  \pst@@divide{#1}{#2}%
  \pst@dimtonum\pst@dimg{#3}%
}
\def\pst@@divide#1#2{%
  \pst@dimg=#1\relax
  \pst@dimh=#2\relax
  \pst@cntg=\pst@dimh
  \pst@cnth=67108863
  \pst@@@divide\pst@@@divide\pst@@@divide\pst@@@divide
  \divide\pst@dimg\pst@cntg%
}
\def\pst@@@divide{%
  \ifnum
    \ifnum\pst@dimg<\z@-\fi\pst@dimg<\pst@cnth
      \multiply\pst@dimg\sixt@@n
    \else
      \divide\pst@cntg\sixt@@n
    \fi%
}
\def\pst@configerr#1{%
  \@pstrickserr{\string#1 not defined in pstricks.con}\@ehpc}
\def\pstVerb#1{\pst@configerr\pstVerb}
\def\pstverb#1{\pst@configerr\pstverb}
\def\pstverbscale{\pst@configerr\pstverbscale}
\def\pstrotate{\pst@configerr\pstrotate}
\def\pstheader#1{\pst@configerr\pstheader}
\def\pstdriver{\pst@configerr\pstdriver}
\@ifundefined{pstcustomize}%
{\def\pstcustomize{ \let\pstcustomize\relax}}{}
\input pstricks.con
\newif\ifPSTricks
\PSTrickstrue

\@ifundefined{pst@def}{\def\pst@def#1<#2>{\@namedef{tx@#1}{#2 }}}{}
\@ifundefined{pst@ATH}{\def\pst@ATH<#1>{}}{}
\pstheader{pstricks.pro}
\def\pst@dict{tx@Dict begin }
\def\pst@theheaders{pstricks.pro}
\def\pst@Verb#1{\pstVerb{\pst@dict #1 end}}
\def\tx@Atan{Atan }
\def\tx@Div{Div }
\def\tx@NET{NET }
\def\tx@Pyth{Pyth }
\def\tx@PtoC{PtoC }
\def\tx@PathLength@{PathLength@ }
\def\tx@PathLength{PathLength }
\pst@dimg=\pstunit\relax
\ifdim\pst@dimg=1bp
\def\pst@stp{.996264 dup scale}
\else
\edef\pst@stp{1 \pst@@dimtonum\pst@dimg\space div dup scale}
\fi
\def\tx@STP{STP }
\def\tx@STV{STV }
\def\pst@number#1{\pst@@dimtonum#1\space}
\def\pst@checknum#1#2{%
  \edef\next{#1}%
  \ifx\next\@empty
    \let\pst@num\z@
  \else
    \expandafter\pst@@checknum\next..\@nil
  \fi
  \ifnum\pst@num=\z@
    \@pstrickserr{Bad number: `#1'. 0 substituted.}\@ehpa
    \def#2{0 }%
  \else
    \edef#2{\ifnum\pst@num=\tw@-\fi\the\pst@cntg.%
    \expandafter\@gobble\the\pst@cnth\space}%
  \fi}
\def\pst@@checknum{%
  \@ifnextchar-%
  {\let\pst@num\tw@\expandafter\pst@@@checknum\@gobble}%
  {\let\pst@num\@ne\pst@@@checknum}%
}
\def\pst@@@checknum#1.#2.#3\@nil{%
\afterassignment\pst@@@@checknum\pst@cntg=0#1\relax\@nil
\afterassignment\pst@@@@checknum\pst@cnth=1#2\relax\@nil}
\def\pst@@@@checknum#1\relax\@nil{%
\ifx\@nil#1\@nil\else\let\pst@num\z@\fi}
\def\pst@getnumii#1 #2 #3\@nil{%
\pst@checknum{#1}\pst@tempg
\pst@checknum{#2}\pst@temph}
\def\pst@getnumiii#1 #2 #3 #4\@nil{%
\pst@checknum{#1}\pst@tempg
\pst@checknum{#2}\pst@temph
\pst@checknum{#3}\pst@tempi}
\def\pst@getnumiv#1 #2 #3 #4 #5\@nil{%
\pst@checknum{#1}\pst@tempg
\pst@checknum{#2}\pst@temph
\pst@checknum{#3}\pst@tempi
\pst@checknum{#4}\pst@tempj}
\def\pst@getdimnum#1 #2 #3\@nil{%
  \pssetlength\pst@dimg{#1}%
  \pst@checknum{#2}\pst@tempg%
}
\def\pst@getscale#1#2{%
  \edef\pst@tempg{#1}%
  \ifx\pst@tempg\@none
    \def#2{}%
  \else
    \pst@expandafter\pst@getnumii{#1 #1} {} {} {}\@nil
    \ifdim\pst@tempg\p@=\z@
      \@pstrickserr{Bad scaling argument `#1'}\@ehpa
      \def#2{}%
    \else
      \ifdim\pst@temph\p@=\z@
        \@pstrickserr{Bad scaling argument `#1'}\@ehpa
        \def#2{}%
      \else
        \edef#2{\pst@tempg\space \pst@temph\space scale }%
      \fi
    \fi
  \fi%
}
\def\pst@getint#1#2{%
  \pst@cntg=#1\relax
  \edef#2{\the\pst@cntg\space}%
}
\begingroup
\catcode`\{=12
\catcode`\}=12
\catcode`\[=1
\catcode`\]=2
\gdef\pslbrace[{ ]
\gdef\psrbrace[} ]
\endgroup
\def\@newcolor#1#2{%
\expandafter\edef\csname #1\endcsname{\noexpand\pst@color{#2}}%
\expandafter\edef\csname\string\color@#1\endcsname{#2}
\ignorespaces}
\def\pst@color#1{%
\def\pst@currentcolor{#1}\pstVerb{#1}\aftergroup\pst@endcolor}
\def\pst@endcolor{\pstVerb{\pst@currentcolor}}
\def\pst@currentcolor{0 setgray}
\def\altcolormode{%
\def\pst@color##1{%
\pstVerb{gsave ##1}\aftergroup\pst@endcolor}%
\def\pst@endcolor{\pstVerb{\pst@grestore}}}
\def\pst@grestore{%
  currentpoint
  matrix currentmatrix
  currentfont
  grestore
  setfont
  setmatrix
  moveto
}
\def\pst@usecolor#1{\csname\string\color@#1\endcsname\space}
\def\newgray#1#2{%
  \pst@checknum{#2}\pst@tempg
  \@newcolor{#1}{\pst@tempg setgray}%
}
\def\newrgbcolor#1#2{%
  \pst@expandafter\pst@getnumiii{#2} {} {} {} {}\@nil
  \@newcolor{#1}{\pst@tempg \pst@temph \pst@tempi setrgbcolor}%
}
\def\newhsbcolor#1#2{%
  \pst@expandafter\pst@getnumiii{#2} {} {} {} {}\@nil
  \@newcolor{#1}{\pst@tempg \pst@temph \pst@tempi sethsbcolor}%
}
\def\newcmykcolor#1#2{%
  \pst@expandafter\pst@getnumiv{#2} {} {} {} {} {}\@nil
  \@newcolor{#1}{\pst@tempg \pst@temph \pst@tempi \pst@tempj setcmykcolor}%
}
\newgray{black}{0}
\newgray{darkgray}{.25}
\newgray{gray}{.5}
\newgray{lightgray}{.75}
\newgray{white}{1}
\newrgbcolor{red}{1 0 0}
\newrgbcolor{green}{0 1 0}
\newrgbcolor{blue}{0 0 1}
\newrgbcolor{yellow}{1 1 0}
\newrgbcolor{cyan}{0 1 1}
\newrgbcolor{magenta}{1 0 1}
\def\psset#1{\@psset#1,\@nil\ignorespaces}
\def\@psset#1,{%
  \@@psset#1==\@nil%
  \@ifnextchar\@nil{\@gobble}{\@psset}%
}
\def\@@psset#1=#2=#3\@nil{%
  \@ifundefined{psset@#1}%
    {\@pstrickserr{Graphics parameter `#1' not defined.}\@ehpa}%
    {\@nameuse{psset@#1}{#2}}%
}%
\def\psset@style#1{%
  \@ifundefined{pscs@#1}%
    {\@pstrickserr{Custom style `#1' undefined}\@ehpa}%
    {\@nameuse{pscs@#1}}%
}
\def\newpsstyle#1#2{\@namedef{pscs@#1}{%
  \def\pst@tempa{#2}%
  \ifx\pst@tempa\@empty\else\psset{#2}\fi}}
\def\@none{none}
\def\pst@getcolor#1#2{%
  \@ifundefined{\string\color@#1}%
    {\@pstrickserr{Color `#1' not defined}\@eha}%
    {\edef#2{#1}}%
}
\newdimen\psunit \psunit 1cm
\newdimen\psxunit \psxunit 1cm
\newdimen\psyunit \psyunit 1cm
\let\psrunit\psunit
\def\pstunit@off{\let\@psunit\ignorespaces\ignorespaces}
\def\pssetlength#1#2{%
  \let\@psunit\psunit
  \afterassignment\pstunit@off
  #1 #2\@psunit%
}
\def\psaddtolength#1#2{%
  \let\@psunit\psunit
  \afterassignment\pstunit@off
  \advance#1 #2\@psunit%
}
\def\pssetxlength#1#2{%
  \let\@psunit\psxunit
  \afterassignment\pstunit@off
  #1 #2\@psunit%
}
\def\pssetylength#1#2{%
  \let\@psunit\psyunit
  \afterassignment\pstunit@off
  #1 #2\@psunit%
}
\def\psset@unit#1{%
  \pssetlength\psunit{#1}%
  \psxunit=\psunit
  \psyunit=\psunit%
}
\def\psset@runit#1{\pssetlength\psrunit{#1}}
\def\psset@xunit#1{\pssetxlength\psxunit{#1}}
\def\psset@yunit#1{\pssetylength\psyunit{#1}}
\def\psset@PstDebug#1{\pst@getint{#1}{\Pst@Debug}}
\psset{PstDebug=0}
\def\pst@getlength#1#2{%
  \pssetlength\pst@dimg{#1}%
  \edef#2{\pst@number\pst@dimg}%
}
\def\pst@@getlength#1#2{%
  \pssetlength\pst@dimg{#1}%
  \edef#2{\number\pst@dimg sp}%
}
\def\pst@getcoor#1#2{\pst@@getcoor{#1}\let#2\pst@coor}
\def\pst@coor{0 0 }
\def\pst@getcoors#1#2{%
  \def\pst@aftercoors{\addto@pscode{#1 \pst@coors }#2}%
  \def\pst@coors{}%
  \pst@@getcoors%
}
\def\pst@@getcoors(#1){%
  \pst@@getcoor{#1}%
  \edef\pst@coors{\pst@coor\pst@coors}%
  \@ifnextchar({\pst@@getcoors}{\pst@aftercoors}%
}
\def\pst@getangle#1#2{\pst@@getangle{#1}\let#2\pst@angle}
\def\pst@angle{0 }
\def\cartesian@coor#1,#2,#3\@nil{%
  \pssetxlength\pst@dimg{#1}%
  \pssetylength\pst@dimh{#2}%
  \edef\pst@coor{\pst@number\pst@dimg \pst@number\pst@dimh}%
}
\def\NormalCoor{%
  \def\pst@@getcoor##1{\pst@expandafter\cartesian@coor{##1},\relax,\@nil}%
  \def\pst@@getangle##1{%
    \pst@checknum{##1}\pst@angle
    \edef\pst@angle{\pst@angle \pst@angleunit}%
  }%
  \def\psput@##1{\pst@@getcoor{##1}\leavevmode\psput@cartesian}%
}
\NormalCoor
\def\degrees{\@ifnextchar[{\@degrees}{\def\pst@angleunit{}}}
\def\@degrees[#1]{%
\pst@checknum{#1}\pst@tempg
\edef\pst@angleunit{360 \pst@tempg div mul }%
\ignorespaces}
\def\radians{\def\pst@angleunit{57.2956 mul }}
\def\pst@angleunit{}
\def\SpecialCoor{%
  \def\pst@@getcoor##1{%
    \begingroup
      \pst@activecoor
      \xdef\pst@tempg{##1}%
    \endgroup
    \expandafter\special@coor\pst@tempg||\@nil%
  }%
  \def\pst@@getangle##1{%
    \begingroup
      \pst@activecoor
      \xdef\pst@tempg{##1}%
    \endgroup
    \expandafter\special@angle\pst@tempg\@empty)\@nil%
  }%
  \def\psput@##1{\pst@@getcoor{##1}\leavevmode\psput@special}%
}
\begingroup
\catcode`\|=13
\catcode`\;=13
\catcode`\!=13
\gdef\pst@activecoor{%
  \def|{\string|}%
  \def;{\string;}%
  \def!{\string!}%
}
\endgroup
\def\special@coor#1|#2|#3\@nil{%
  \ifx#3|\relax
    \mixed@coor{#1}{#2}%
  \else
    \special@@coor#1;;\@nil
  \fi%
}
\def\special@@coor#1{%
  \ifcat#1a\relax
    \def\next{\node@coor#1}%
  \else
    \ifx#1[\relax
      \def\next{\Node@coor[}%
    \else
      \ifx#1!\relax
        \def\next{\raw@coor}%
      \else
        \def\next{\special@@@coor#1}%
      \fi
    \fi
  \fi
  \next%
}
\def\special@@@coor#1;#2;#3\@nil{%
  \ifx#3;\relax
    \polar@coor{#1}{#2}%
  \else
    \cartesian@coor#1,\relax,\@nil
  \fi%
}
\def\mixed@coor#1#2{%
  \begingroup
    \special@@coor#1;;\@nil
    \let\pst@tempa\pst@coor
    \special@@coor#2;;\@nil
    \xdef\pst@tempg{\pst@tempa pop \pst@coor exch pop }%
  \endgroup
  \let\pst@coor\pst@tempg%
}
\def\polar@coor#1#2{%
  \pssetlength\pst@dimg{#1}%
  \pst@@getangle{#2}%
  \edef\pst@coor{\pst@number\pst@dimg \pst@angle \tx@PtoC}%
}
\def\raw@coor#1;#2\@nil{%
  \edef\pst@coor{%
    #1 \pst@number\psyunit mul exch \pst@number\psxunit mul exch }%
}
\def\node@coor#1\@nil{%
  \@pstrickserr{You must load `pst-node.tex' to use node coordinates.}\@ehps
  \def\pst@coor{0 0 }%
}
\def\Node@coor{\node@coor}
\def\special@angle#1#2)#3\@nil{%
\ifx!#1\relax
\edef\pst@angle{#2 \pst@angleunit}%
\else
\ifx(#1\relax
\pst@@getcoor{#2}%
\edef\pst@angle{\pst@coor exch \tx@Atan}%
\else
\pst@checknum{#1#2}\pst@angle
\edef\pst@angle{\pst@angle \pst@angleunit}%
\fi
\fi}
\def\Cartesian{%
  \def\cartesian@coor##1,##2,##3\@nil{%
    \pssetxlength\pst@dimg{##1}%
    \pssetylength\pst@dimh{##2}%
    \edef\pst@coor{\pst@number\pst@dimg \pst@number\pst@dimh}%
  }%
  \@ifnextchar({\Cartesian@}{}%
}
\def\Cartesian@(#1,#2){%
  \pssetxlength\psxunit{#1}%
  \pssetylength\psyunit{#2}%
  \ignorespaces%
}
\def\Polar{%
  \def\psput@cartesian{\psput@special}%
  \def\cartesian@coor##1,##2,##3\@nil{\polar@coor{##1}{##2}}%
}%
\def\psset@origin#1{%
  \pst@@getcoor{#1}%
  \edef\psk@origin{\pst@coor T }}
\def\psk@origin{}
\newif\ifswapaxes
\def\psset@swapaxes#1{%
  \@nameuse{@pst#1}%
  \if@pst\def\psk@swapaxes{-90 rotate -1 1 scale }%
  \else\def\psk@swapaxes{}%
  \fi%
}
\psset@swapaxes{false}
\newif\ifshowpoints
\def\psset@showpoints#1{\@nameuse{showpoints#1}}
\psset@showpoints{false}
\let\pst@setrepeatarrowsflag\relax
\def\psset@border#1{%
\pst@getlength{#1}\psk@border
\pst@setrepeatarrowsflag}
\psset@border{0pt}
\def\psset@bordercolor#1{\pst@getcolor{#1}\psbordercolor}
\psset@bordercolor{white}
\newif\ifpsdoubleline
\def\psset@doubleline#1{%
  \@nameuse{psdoubleline#1}%
  \pst@setrepeatarrowsflag}
\psset@doubleline{false}
\def\psset@doublesep#1{\def\psdoublesep{#1}}
\psset@doublesep{1.25\pslinewidth}
\def\psset@doublecolor#1{\pst@getcolor{#1}\psdoublecolor}
\psset@doublecolor{white}
\newif\ifpsshadow
\def\psset@shadow#1{%
  \@nameuse{psshadow#1}%
  \pst@setrepeatarrowsflag}
\psset@shadow{false}
\def\psset@shadowsize#1{\pst@getlength{#1}\psk@shadowsize}
\psset@shadowsize{3pt}
\def\psset@shadowangle#1{\pst@getangle{#1}\psk@shadowangle}
\psset@shadowangle{-45}
\def\psset@shadowcolor#1{\pst@getcolor{#1}\psshadowcolor}
\psset@shadowcolor{darkgray}
\def\pst@repeatarrowsflag{\z@}
\def\pst@setrepeatarrowsflag{%
  \edef\pst@repeatarrowsflag{%
    \ifdim\psk@border\p@>\z@ 1\else\ifpsdoubleline 1\else
      \ifpsshadow 1\else \z@\fi\fi\fi}}
\def\psls@none{}
\newdimen\pslinewidth
\def\psset@linewidth#1{\pssetlength\pslinewidth{#1}}
\psset@linewidth{.8pt}
\def\psset@linecolor#1{\pst@getcolor{#1}\pslinecolor}
\psset@linecolor{black}
\def\psls@solid{0 setlinecap stroke }
\def\pst@missing{%
  \z@
  \@pstrickserr{Missing number or dimension. 0 substituted}\@ehpa}
%
\def\pst@empty{\z@}
\def\psset@dash#1{
  \pst@expandafter\psset@@dash{#1} {\pst@empty} {\pst@empty} %
                  {\pst@missing} {\pst@missing} {}\@nil
  \edef\psk@dash{\pst@number\pst@dimg \pst@number\pst@dimh 
                 \pst@number\pst@dimc \pst@number\pst@dimd}%
}
\def\psset@@dash#1 #2 #3 #4 #5\@nil{%
  \pssetlength\pst@dimg{#1}%
  \pssetlength\pst@dimh{#2}%
  \pssetlength\pst@dimc{#3}%
  \pssetlength\pst@dimd{#4}%
}
\psset@dash{5pt 3pt 0pt 0pt}
\newif\ifpsdashadjust
\def\psset@dashadjust#1{\@nameuse{psdashadjust#1}}
\psset@dashadjust{true}
\def\tx@DashLine{DashLine }
\def\psls@dashed{%
  \psk@linecap\space setlinecap
  \ifpsdashadjust
    \psk@dash \@ifundefined{pst@linetype}{2}{\pst@linetype}\space \tx@DashLine
  \else
    [ \psk@dash ] 0 setdash stroke
  \fi}
\def\psset@dotsep#1{\pst@getlength{#1}\psk@dotsep}
\psset@dotsep{3pt}
\def\tx@DotLine{DotLine }
\def\psls@dotted{%
  \ifpsdashadjust
    \psk@dotsep \pst@linetype\space \tx@DotLine
  \else
    [ 0 \psk@dotsep CLW add ] 0 setdash 1 setlinecap stroke
  \fi%
}
\def\psset@linestyle#1{%
  \@ifundefined{psls@#1}%
    {\@pstrickserr{Line style `#1' not defined}\@eha}%
    {\edef\pslinestyle{#1}}%
}
\psset@linestyle{solid}
\def\psset@linecap#1{%
  \def\psk@linecap{0}%
  \ifnum#1>-1
    \ifnum#1<3
      \pst@getint{#1}\psk@linecap%
  \fi\fi%
}
\psset{linecap=0}

\def\psfs@none{}
\def\psset@fillcolor#1{\pst@getcolor{#1}\psfillcolor}
\psset@fillcolor{white}
\def\psfs@solid{\pst@fill{\pst@usecolor\psfillcolor fill}}
\def\psfs@eofill{\pst@fill{\pst@usecolor\psfillcolor eofill}}
\def\psset@hatchwidth#1{\pst@getlength{#1}\psk@hatchwidth}
\psset@hatchwidth{.8pt}
\def\psset@hatchsep#1{\pst@getlength{#1}\psk@hatchsep}
\psset@hatchsep{4pt}
\def\psset@hatchcolor#1{\pst@getcolor{#1}\pshatchcolor}
\psset@hatchcolor{black}
\def\psset@hatchangle#1{\pst@getangle{#1}\psk@hatchangle}
\psset@hatchangle{45}
\def\psset@hatchsepinc#1{\pst@getlength{#1}\psk@hatchsepinc}
\def\psset@hatchwidthinc#1{\pst@getlength{#1}\psk@hatchwidthinc}
\psset@hatchwidthinc{0pt}
\psset@hatchsepinc{0pt}
\def\pst@linefill{%
  \psk@hatchangle rotate
  \psk@hatchwidth SLW
  \pst@usecolor\pshatchcolor
  \psk@hatchsep 
  \psk@hatchsepinc
  \psk@hatchwidthinc
  \tx@LineFill }
%
\def\psfs@vlines{\pst@fill\pst@linefill}
\@namedef{psfs@vlines*}{\psfs@solid \psfs@vlines}
\def\psfs@hlines{\pst@fill{90 rotate \pst@linefill}}
\@namedef{psfs@hlines*}{\psfs@solid \psfs@hlines}
\def\psfs@crosshatch{\psfs@vlines \psfs@hlines}
\@namedef{psfs@crosshatch*}{\psfs@solid \psfs@vlines \psfs@hlines}
\def\tx@LineFill{LineFill }
\def\psset@fillstyle#1{%
  \edef\pst@tempg{#1}\def\pst@temph{none}%
  \ifx\pst@tempg\pst@temph
     \let\psk@fillstyle\relax
  \else
    \@ifundefined{psfs@#1}%
    {\@pstrickserr{Undefined fill style: `#1'}\@eha}%
    {\edef\psk@fillstyle{\expandafter\noexpand\csname psfs@#1\endcsname}}%
  \fi%
}
\def\psset@addfillstyle#1{%
  \@ifundefined{psfs@#1}%
    {\@pstrickserr{Undefined fill style: `#1'}\@eha}%
    {\edef\psk@fillstyle{%
      \expandafter\noexpand\psk@fillstyle
      \expandafter\noexpand\csname psfs@#1\endcsname}%
    }%
}
\psset@fillstyle{none}
\def\psset@arrows#1{%
  \begingroup
    \pst@activearrows
    \xdef\pst@tempg{#1}%
  \endgroup
  \expandafter\psset@@arrows\pst@tempg\@empty-\@empty\@nil
  \if@pst\else
    \@pstrickserr{Bad arrows specification: #1}\@ehpa
  \fi%
}
\def\psset@@arrows#1-#2\@empty#3\@nil{%
  \@psttrue
  \def\next##1,#1-##2,##3\@nil{\def\pst@tempg{##2}}%
  \expandafter\next\pst@arrowtable,#1-#1,\@nil
  \@ifundefined{psas@\pst@tempg}%
    {\@pstfalse\def\psk@arrowA{}}%
    {\let\psk@arrowA\pst@tempg}%
  \@ifundefined{psas@#2}%
    {\@pstfalse\def\psk@arrowB{}}%
    {\def\psk@arrowB{#2}}%
}
\def\psk@arrowA{}
\def\psk@arrowB{}
\def\pst@arrowtable{,<->,<<->>,>-<,>>-<<,(-),[-],)-(,]-[,|>-<|} 
\begingroup
  \catcode`\<=13
  \catcode`\>=13
  \catcode`\|=13
  \gdef\pst@activearrows{\def<{\string<}\def>{\string>}\def|{\string|}}
\endgroup
\def\tx@BeginArrow{BeginArrow }
\def\tx@EndArrow{EndArrow }
\def\psset@arrowscale#1{
  \pst@@arrowscale@i#1 \@nil
  \pst@getscale{\pst@arrowscale}\psk@arrowscale}
\def\pst@@arrowscale@i#1 #2\@nil{\edef\pst@arrowscale{#1}}
\psset@arrowscale{1}
\def\psset@arrowsize#1{%
  \pst@expandafter\pst@getdimnum{#1} 0 {} {}\@nil
  \edef\psk@arrowsize{\pst@number\pst@dimg \pst@tempg}%
}
\psset@arrowsize{1.5pt 2}
\def\psset@arrowlength#1{\pst@checknum{#1}\psk@arrowlength}
\psset@arrowlength{1.4}
\def\psset@arrowinset#1{\pst@checknum{#1}\psk@arrowinset}%
\psset@arrowinset{.4}
\def\tx@Arrow{Arrow }
\@namedef{psas@<|}{%
    \psk@tbarsize\space \tx@Tbar
    0 CLW 2 div T
    newpath
    true \psk@arrowinset\space \psk@arrowlength\space \psk@arrowsize\space \tx@Arrow%
}
\def\tx@BracketOut{BracketOut }
\@namedef{psas@[}{%
  /BracketOut {%
  CLW mul add dup CLW sub 2 div
  /x ED mul neg
  /y ED
  /z CLW 2 div def
  x neg y moveto
  x neg CLW 2 div L x CLW 2 div L x y L stroke 0 CLW moveto } def
  \psk@bracketlength\space \psk@tbarsize\space \tx@BracketOut
}
\def\tx@RoundBracketOut{RoundBracketOut }
\@namedef{psas@(}{%
  /RoundBracketOut {%
    CLW mul add dup 2 div
    /x ED mul neg
    /y ED
    /mtrx CM def
    0 CLW
    2 div T x y mul 0 ne { x y scale } if
    1 1 moveto
    .85 .5 .35 0 0 0 curveto
    -.35 0 -.85 .5 -1 1 curveto
    mtrx setmatrix stroke 0 CLW moveto } def
  \psk@rbracketlength\space \psk@tbarsize\space \tx@RoundBracketOut
}
\@namedef{psas@>}{%
  false \psk@arrowinset \psk@arrowlength \psk@arrowsize \tx@Arrow
}
\@namedef{psas@>>}{%
  false \psk@arrowinset \psk@arrowlength \psk@arrowsize \tx@Arrow
  0 h T
  gsave
  newpath
  false \psk@arrowinset \psk@arrowlength \psk@arrowsize \tx@Arrow
  CP
  grestore
  CP newpath moveto
  2 copy
  L
  stroke
  moveto
}
\@namedef{psas@<}{true \psk@arrowinset \psk@arrowlength \psk@arrowsize \tx@Arrow}
\@namedef{psas@<<}{%
  true \psk@arrowinset \psk@arrowlength \psk@arrowsize \tx@Arrow
  CP newpath moveto 0 a neg L stroke 0 h neg T
  false \psk@arrowinset \psk@arrowlength \psk@arrowsize \tx@Arrow
}
\def\psset@tbarsize#1{%
  \pst@expandafter\pst@getdimnum{#1} 0 {} {}\@nil
  \edef\psk@tbarsize{\pst@number\pst@dimg \pst@tempg}%
}
\psset@tbarsize{2pt 5}
\def\tx@Tbar{Tbar }
\@namedef{psas@|}{\psk@tbarsize \tx@Tbar}
\@namedef{psas@|*}{0 CLW -2 div T \psk@tbarsize \tx@Tbar}
\@namedef{psas@>|}{%
  \psk@tbarsize \tx@Tbar
  0 CLW 2 div T
  newpath
  false \psk@arrowinset \psk@arrowlength \psk@arrowsize \tx@Arrow
}
\@namedef{psas@>|*}{%
  0 CLW -2 div T
  \psk@tbarsize \tx@Tbar
  0 CLW 2 div T
  newpath
  false \psk@arrowinset \psk@arrowlength \psk@arrowsize \tx@Arrow
}
\edef\pst@arrowtable{\pst@arrowtable,|<*->|*,|<->|}
\def\psset@bracketlength#1{\pst@checknum{#1}\psk@bracketlength}
\psset@bracketlength{.15}
\def\tx@Bracket{Bracket }
\@namedef{psas@]}{\psk@bracketlength \psk@tbarsize \tx@Bracket}
\def\psset@rbracketlength#1{\pst@checknum{#1}\psk@rbracketlength}
\psset@rbracketlength{.15}
\def\tx@RoundBracket{RoundBracket }
\@namedef{psas@)}{\psk@rbracketlength \psk@tbarsize \tx@RoundBracket}
\def\psas@c{1 \psas@@c}
\def\psas@cc{0 CLW 2 div T 1 \psas@@c}
\def\psas@C{2 \psas@@c}
\def\psas@@c{%
  setlinecap
  0 0 moveto
  0 0.1 L 
  stroke
  0 0 moveto
}
\def\psas@{}
\psset@arrows{-}
\def\tx@SD{SD }
\def\tx@EndDot{EndDot }
\def\psas@oo{{\pst@usecolor\psfillcolor true} true \psk@dotsize \tx@EndDot}
\def\psas@o{{\pst@usecolor\psfillcolor true} false \psk@dotsize \tx@EndDot}
\@namedef{psas@**}{{false} true \psk@dotsize \tx@EndDot}
\@namedef{psas@*}{{false} false \psk@dotsize \tx@EndDot}
\def\pst@par{}
\def\addto@par#1{%
\ifx\pst@par\@empty
\def\pst@par{#1}%
\else
\expandafter\def\expandafter\pst@par\expandafter{\pst@par,#1}%
\fi}
\def\addbefore@par#1{%
\ifx\pst@par\@empty
\def\pst@par{#1}%
\else
\toks@{#1}%
\pst@toks\expandafter{\pst@par}%
\edef\pst@par{\the\toks@,\the\pst@toks}%
\fi}
\def\use@par{%
  \ifx\pst@par\@empty\else
    \expandafter\@psset\pst@par,\@nil
    \def\pst@par{}%
  \fi%
}
\def\pst@object#1{%
  \pst@ifstar{%
    \@ifnextchar[{\pst@@object{#1}}{\def\pst@par{}\@nameuse{#1@i}}}%
}
\def\pst@@object#1[#2]{%
  \def\pst@par{#2}%
  \@ifnextchar+{\@nameuse{#1@i}}{\@nameuse{#1@i}}}
\def\newpsobject#1#2#3{%
\@ifundefined{#2@i}%
{\@pstrickserr{Graphics object `#2' not defined}\@eha}{%
\@namedef{#1}{\pst@object{#1}}%
\@namedef{#1@i}{\addbefore@par{#3}\@nameuse{#2@i}}}%
\ignorespaces}
\def\pst@getarrows#1{\@ifnextchar({#1}{\pst@@getarrows{#1}}}
\def\pst@@getarrows#1#2{%
  \def\pst@tempa{#2}
  \ifx\pst@tempa\@empty\addto@par{arrows=-}\else\addto@par{arrows=#2}\fi#1}
%
\def\begin@ClosedObj{%
  \leavevmode
  \pst@killglue
  \begingroup
    \use@par
    \solid@star
    \ifpsdoubleline \pst@setdoublesep \fi
    \init@pscode%
}
\def\end@ClosedObj{%
  \ifpsshadow \pst@closedshadow \fi
  \ifdim\psk@border\p@>\z@ \pst@addborder \fi
  \psk@fillstyle
  \pst@stroke
  \ifpsdoubleline \pst@doublestroke \fi
  \ifshowpoints
  \pst@OpenShowPoints
  \fi
  \use@pscode
  \endgroup
  \ignorespaces%
}
\def\begin@OpenObj{%
  \begin@ClosedObj
  \let\pst@linetype\pst@arrowtype
  \pst@addarrowdef%
}
\def\begin@AltOpenObj{%
  \begin@ClosedObj
  \def\pst@repeatarrowsflag{\z@}%
  \def\pst@linetype{0}}
\def\end@OpenObj{%
  \ifpsshadow \pst@openshadow \fi
  \ifdim\psk@border\p@>\z@ \pst@addborder \fi
  \psk@fillstyle
  \pst@stroke
  \ifpsdoubleline \pst@doublestroke \fi
  \ifnum\pst@repeatarrowsflag>\z@ \pst@repeatarrows \fi
  \ifshowpoints \pst@OpenShowPoints \fi
  \use@pscode
  \endgroup
  \ignorespaces}
\def\begin@SpecialObj{%
\leavevmode
\pst@killglue
\begingroup
\use@par
\init@pscode}
\def\end@SpecialObj{%
\use@pscode
\endgroup
\ignorespaces}
\def\pst@code{}%
\def\init@pscode{%
  \addto@pscode{%
    \pst@number\pslinewidth SLW
    \pst@usecolor\pslinecolor}%
}
\def\addto@pscode#1{\xdef\pst@code{\pst@code#1\space}}
\def\use@pscode{%
  \pstverb{%
    \pst@dict
    \tx@STP
    \pst@newpath
    \psk@origin
    \psk@swapaxes
    \pst@code
    end
  }%
  \gdef\pst@code{}%
}
\def\pst@newpath{newpath }
\def\pst@@killglue{\unskip\ifdim\lastskip>\z@\expandafter\pst@@killglue\fi}
\def\KillGlue{\let\pst@killglue\pst@@killglue}
\def\DontKillGlue{\let\pst@killglue\relax}
\DontKillGlue
\def\solid@star{%
  \if@star
    \pslinewidth=\z@
    \psdoublelinefalse
    \def\pslinestyle{none}%
    \def\psk@fillstyle{\psfs@solid}%
    \let\psfillcolor\pslinecolor
  \fi}
\def\pst@setdoublesep{%
\pst@getlength\psdoublesep\psdoublesep
\pslinewidth=2\pslinewidth
\advance\pslinewidth\psdoublesep\p@
\let\pst@setdoublesep\relax}
\def\tx@Shadow{Shadow }
\def\pst@closedshadow{%
  \addto@pscode{%
    gsave
    \psk@shadowsize \psk@shadowangle \tx@PtoC
    \tx@Shadow
    \pst@usecolor\psshadowcolor
    gsave fill grestore
    stroke
    grestore
    gsave
    \pst@usecolor\psfillcolor
    gsave fill grestore
    stroke
    grestore}}
\def\pst@openshadow{%
  \addto@pscode{%
    gsave
    \psk@shadowsize \psk@shadowangle \tx@PtoC
    \tx@Shadow
    \pst@usecolor\psshadowcolor
    \ifx\psk@fillstyle\relax\else
      gsave fill grestore
    \fi
    stroke}%
  \pst@repeatarrows%
  \addto@pscode{grestore}%
  \ifx\psk@fillstyle\relax\else
    \addto@pscode{%
      gsave
      \pst@usecolor\psfillcolor
      gsave fill grestore
      stroke
      grestore}%
  \fi}
\def\pst@addborder{%
  \addto@pscode{%
    gsave
    \psk@border 2 mul
    CLW add SLW
    \pst@usecolor\psbordercolor
    stroke
    grestore}}
\def\pst@stroke{%
  \ifx\pslinestyle\@none\else
    \addto@pscode{%
      gsave
      \pst@number\pslinewidth SLW
      \pst@usecolor\pslinecolor
      \@nameuse{psls@\pslinestyle}
      grestore}%
  \fi}
\def\pst@fill#1{\addto@pscode{gsave #1 grestore}}%
\def\pst@doublestroke{%
    \addto@pscode{%
      gsave
      \psdoublesep SLW
      \pst@usecolor\psdoublecolor
      stroke
      grestore
    }}
\def\pst@arrowtype{%
\ifx\psk@arrowB\@empty 0 \else -2 \fi
\ifx\psk@arrowA\@empty 0 \else -1 \fi
add}
\def\pst@addarrowdef{%
\addto@pscode{%
/ArrowA {
\ifx\psk@arrowA\@empty
  \pst@oplineto
\else
  \pst@arrowdef{A}
moveto
\fi
} def
/ArrowB {
\ifx\psk@arrowB\@empty \else \pst@arrowdef{B} \fi
} def}}
\def\pst@arrowdef#1{%
\ifnum\pst@repeatarrowsflag>\z@
/Arrow#1c [ 6 2 roll ] cvx def Arrow#1c
\fi
\tx@BeginArrow
\psk@arrowscale
\@nameuse{psas@\@nameuse{psk@arrow#1}}
\tx@EndArrow}
\def\pst@repeatarrows{%
\addto@pscode{%
gsave
\ifx\psk@arrowA\@empty\else
ArrowAc ArrowA pop pop
\fi
\ifx\psk@arrowB\@empty\else
ArrowBc ArrowB pop pop pop pop
\fi
grestore}}
\def\pst@OpenShowPoints{%
  \addto@pscode{%
    gsave
    \psk@dotsize
    \@nameuse{psds@\psk@dotstyle}
    newpath
    Points aload length 2 div 2 sub cvi /N ED
    N 0 ge
      { \ifx\psk@arrowA\@empty Dot \else pop pop \fi 
        N { Dot } repeat 
	\ifx\psk@arrowB\@empty Dot \else pop pop \fi }
      { N 2 mul { pop } repeat } ifelse
    grestore
}}
\def\pscustom{\pst@object{pscustom}}
\long\def\pscustom@i#1{%
  \begin@SpecialObj
  \solid@star
  \let\pst@ifcustom\iftrue
  \let\begin@ClosedObj\begin@CustomObj
  \let\end@ClosedObj\endgroup
  \def\begin@OpenObj{\begin@CustomObj\pst@addarrowdef}%
  \let\end@OpenObj\endgroup
  \let\begin@AltOpenObj\begin@CustomObj
  \def\begin@SpecialObj{%
    \begingroup
    \pst@misplaced{special graphics object}%
    \def\addto@pscode####1{}
    \let\end@SpecialObj\endgroup}%
    \def\@multips(##1)(##2)##3##4{\pst@misplaced\multips}%
    \def\psclip##1{\pst@misplaced\psclip}%
    \def\pst@repeatarrowsflag{\z@}%
    \let\pst@setrepeatarrowsflag\relax
    \showpointsfalse
    \let\showpointstrue\relax
    \def\pst@linetype{\pslinetype}%
    \let\psset@liftpen\psset@@liftpen
    \psset@liftpen{\z@}%
    \def\pst@cp{/currentpoint load stopped pop }%
    \def\pst@oplineto{/lineto load stopped { moveto } if }%
    \def\pst@optcp##1##2{%
    \ifnum##1=\z@\def##2{/currentpoint load stopped { 0 0 } if }\fi}%
    \let\caddto@pscode\addto@pscode
    \def\cuse@par##1{{\use@par##1}}%
    \the\pst@customdefs
    \setbox\pst@hbox=\hbox{#1}%
    \psk@fillstyle
    \pst@stroke
  \end@SpecialObj}
\def\begin@CustomObj{%
  \begingroup
  \use@par
  \addto@pscode{%
    \pst@number\pslinewidth SLW
    \pst@usecolor\pslinecolor
  }%
}
\def\pst@oplineto{moveto }
\def\pst@cp{}
\def\pst@optcp#1#2{}
\def\psset@liftpen#1{}
\def\psset@@liftpen#1{%
  \ifcase#1\relax
    \def\psk@liftpen{\z@}%
    \def\pst@cp{/currentpoint load stopped pop }%
    \def\pst@oplineto{/lineto load stopped { moveto } if }%
  \or
    \def\psk@liftpen{1}%
    \def\pst@cp{}%
    \def\pst@oplineto{/lineto load stopped { moveto } if }%
  \or
    \def\psk@liftpen{2}%
    \def\pst@cp{}%
    \def\pst@oplineto{moveto }%
  \fi%
}
\psset@liftpen{0}
\def\psk@liftpen{-1}
\def\psset@linetype#1{%
  \pst@getint{#1}\pslinetype
  \ifnum\pst@cntg<-3
    \@pstrickserr{linetype must be greater than -3}\@ehpa
    \def\pslinetype{2}%
  \fi%
}
\psset@linetype{2}
%
\def\caddto@pscode#1{%
    \@pstrickserr{Command can only be used in \string\pscustom}\@ehpa%
}
\let\cuse@par\caddto@pscode
\def\tx@MSave{%
     /msavematrx
         [ tx@Dict /msavematrx known 
             { msavematrx aload pop } if
             CM 
         ]
     def
    msavematrx
}
\def\tx@MRestore{
     tx@Dict /msavematrx known { length 0 gt } { false } ifelse
         { msavematrx aload pop setmatrix } if
}
\newtoks\pst@customdefs
\pst@customdefs{%
  \def\newpath{\addto@pscode{newpath}}%
  \def\moveto(#1){\pst@@getcoor{#1}\addto@pscode{\pst@coor moveto}}%
  \def\closepath{\addto@pscode{closepath}}%
  \def\gsave{\begingroup\addto@pscode{gsave}}%
  \def\grestore{\endgroup\addto@pscode{grestore}}%
  \def\translate(#1){\pst@@getcoor{#1}\addto@pscode{\pst@coor translate}}%
  \def\rotate#1{\pst@@getangle{#1}\addto@pscode{\pst@angle rotate}}%
  \def\scale#1{\pst@getscale{#1}\pst@tempg\addto@pscode{\pst@tempg}}%
  \def\msave{\addto@pscode{\tx@MSave}}%
  \def\mrestore{\addto@pscode{\tx@MRestore}}%
  \def\swapaxes{\addto@pscode{-90 rotate -1 1 scale}}%
  \def\stroke{\pst@object{stroke}}%
  \def\fill{\pst@object{fill}}%
  \def\openshadow{\pst@object{openshadow}}%
  \def\closedshadow{\pst@object{closedshadow}}%
  \def\movepath(#1){\pst@@getcoor{#1}\addto@pscode{\pst@coor \tx@Shadow}}%
  \def\lineto{\pst@onecoor{lineto}}%
  \def\rlineto{\pst@onecoor{rlineto}}%
  \def\curveto{\pst@threecoor{curveto}}%
  \def\rcurveto{\pst@threecoor{rcurveto}}%
  \def\code#1{\addto@pscode{#1}}%
  \def\coor(#1){\pst@@getcoor{#1}\addto@pscode\pst@coor\@ifnextchar({\coor}{}}%
  \def\rcoor{\pst@getcoors{}{}}%
  \def\dim#1{\pssetlength\pst@dimg{#1}\addto@pscode{\pst@number\pst@dimg}}%
  \def\setcolor#1{%
    \@ifundefined{\string\color@#1}{}{\addto@pscode{\pst@usecolor{#1}}}}
  \def\arrows#1{{\psset@arrows{#1}\pst@addarrowdef}}%
  \let\file\pst@rawfile
} 
\def\closedshadow@i{\cuse@par\pst@closedshadow}
\def\openshadow@i{\cuse@par\pst@openshadow}
\def\stroke@i{\cuse@par\pst@stroke}%
\def\fill@i{\cuse@par\psk@fillstyle}%
\def\pst@onecoor#1(#2){%
\pst@@getcoor{#2}%
\addto@pscode{\pst@coor #1}}
\def\pst@threecoor#1(#2)#3(#4)#5(#6){%
\begingroup
\pst@getcoor{#2}\pst@tempa
\pst@getcoor{#4}\pst@tempb
\pst@getcoor{#6}\pst@tempc
\addto@pscode{\pst@tempa \pst@tempb \pst@tempc #1}%
\endgroup}
\def\pst@rawfile#1{%
\begingroup
\def\do##1{\catcode`##1=12\relax}"
\dospecials
\catcode`\%=14
\pst@@rawfile{#1}%
\endgroup}
\def\pst@@rawfile#1{%
\immediate\openin1 #1
\ifeof1
\@pstrickserr{File `#1' not found}\@ehpa
\else
\immediate\read1 to \pst@tempg
\loop
\ifeof1 \@pstfalse\else\@psttrue\fi
\if@pst
\addto@pscode\pst@tempg
\immediate\read1 to \pst@tempg
\repeat
\fi
\immediate\closein1\relax}
\def\tx@NArray{NArray }
\def\tx@NArray{NArray }
\def\tx@Line{Line }
\def\tx@Arcto{Arcto }
\def\tx@CheckClosed{CheckClosed }
\def\tx@Polygon{Polygon }
\def\psset@gangle#1{\pst@getangle{#1}\psk@gangle}
\psset@gangle{0}
\def\tx@Diamond{Diamond }
\def\psdiamond{\pst@object{psdiamond}}
\def\psdiamond@i(#1){%
\@ifnextchar({\psdiamond@ii(#1)}{\psdiamond@ii(0,0)(#1)}}
\def\psdiamond@ii(#1)(#2){%
  \begin@ClosedObj
  \pst@getcoor{#1}\pst@tempa
  \pst@getcoor{#2}\pst@tempb
  \addto@pscode{%
    \psline@iii
    pop
    \psk@dimen
    \pst@tempb
    \psk@gangle
    \pst@tempa
    \tx@Diamond
  }%
  \def\pst@linetype{4}%
  \end@ClosedObj}
\def\tx@Triangle{Triangle }
\def\pstriangle{\pst@object{pstriangle}}
\def\pstriangle@i(#1){%
  \@ifnextchar({\pstriangle@ii(#1)}{\pstriangle@ii(0,0)(#1)}}
\def\pstriangle@ii(#1)(#2){%
  \begin@ClosedObj
  \pst@getcoor{#1}\pst@tempa
  \pst@getcoor{#2}\pst@tempb
  \addto@pscode{%
    \psline@iii
    pop			    
    \psk@dimen		    
    \pst@tempb
    \psk@gangle		    
    \pst@tempa
    \tx@Triangle
  }%
  \def\pst@linetype{2}%
  \end@ClosedObj}
\def\tx@CCA{CCA }
\def\tx@CCA{CCA }
\def\tx@CC{CC }
\def\tx@IC{IC }
\def\tx@BOC{BOC }
\def\tx@NC{NC }
\def\tx@EOC{EOC }
\def\tx@BAC{BAC }
\def\tx@NAC{NAC }
\def\tx@EAC{EAC }
\def\tx@OpenCurve{OpenCurve }
\def\tx@AltCurve{AltCurve }
\def\tx@ClosedCurve{ClosedCurve }
\def\psset@curvature#1{%
\edef\pst@tempg{#1 }%
\expandafter\psset@@curvature\pst@tempg * * * \@nil}
\def\psset@@curvature#1 #2 #3 #4\@nil{%
\pst@checknum{#1}\pst@tempg
\pst@checknum{#2}\pst@temph
\pst@checknum{#3}\pst@tempi
\edef\psk@curvature{\pst@tempg \pst@temph \pst@tempi}}
\psset@curvature{1 .1 0}
\def\pscurve{\pst@object{pscurve}}
\def\pscurve@i{%
  \pst@getarrows{%
    \begin@OpenObj
      \pst@getcoors[\pscurve@ii%
    }%
}
\def\pscurve@ii{%
  \addto@pscode{%
    \pst@cp
    \psk@curvature\space /c ED /b ED /a ED
    \ifshowpoints true \else false \fi
    \tx@OpenCurve%
  }%
  \end@OpenObj%
}
\def\psecurve{\pst@object{psecurve}}
\def\psecurve@i{%
\pst@getarrows{%
\begin@OpenObj
\pst@getcoors[\psecurve@ii}}
\def\psecurve@ii{%
\addto@pscode{%
\psk@curvature\space /c ED /b ED /a ED
\ifshowpoints true \else false \fi
\tx@AltCurve}%
\end@OpenObj}
\def\psccurve{\pst@object{psccurve}}
\def\psccurve@i{%
\begin@ClosedObj
\pst@getcoors[\psccurve@ii}
\def\psccurve@ii{%
\addto@pscode{%
\psk@curvature\space /c ED /b ED /a ED
\ifshowpoints true \else false \fi
\tx@ClosedCurve}%
\def\pst@linetype{1}%
\end@ClosedObj}
\def\psset@dotsize#1{%
\pst@expandafter\pst@getdimnum{#1} 0 {} {}\@nil
\edef\psk@@dotsize{\pst@number\pst@dimg}%
\let\psk@@@dotsize\pst@tempg
\edef\psk@dotsize{%
/DS \psk@@dotsize \psk@@@dotsize CLW mul add 2 div def }}
\psset@dotsize{2pt 2}
\def\psset@dotscale#1{%
\pst@getscale{#1}\psk@dotscale
\ifx\psk@dotscale\@empty
\def\psk@xdotscale{1 }%
\def\psk@ydotscale{1 }%
\else
\let\psk@xdotscale\pst@tempg
\let\psk@ydotscale\pst@temph
\fi}
\def\pst@Getangle#1#2{%
\pst@getangle{#1}\pst@tempg
\def\pst@temph{0. }%
\ifx\pst@tempg\pst@temph
\def#2{}%
\else
\edef#2{\pst@tempg\space rotate }%
\fi}
\def\psset@dotangle#1{%
\pst@getangle{#1}\psk@@dotangle
\ifdim\psk@@dotangle\p@=\z@
\let\psk@dotangle\@empty
\else
\edef\psk@dotangle{\psk@@dotangle rotate }%
\fi}
\psset@dotangle{0}
\def\pst@getdotsize{%
\pst@dimg=\psk@@@dotsize\pslinewidth
\advance\pst@dimg\psk@@dotsize\p@
\pst@dimh=\psk@ydotscale\pst@dimg
\pst@dimg=\psk@xdotscale\pst@dimg
\divide\pst@dimh 2
\divide\pst@dimg 2\relax}
\psset@dotscale{1}
\def\psdot{\pst@object{psdot}}
\def\psdot@i{\@ifnextchar({\psdot@ii}{\psdot@ii(\z@,\z@)}}
\def\psdot@ii(#1){%
  \begin@SpecialObj
  \solid@star
  \pst@@getcoor{#1}%
  \addto@pscode{%
    \psk@dotsize
    \@nameuse{psds@\psk@dotstyle}%
    \pst@coor Dot}%
  \end@SpecialObj}
\def\psdots{\pst@object{psdots}}
\def\psdots@i{%
\begin@SpecialObj
\pst@getcoors[\psdots@ii}
\def\psdots@ii{%
\addto@pscode{false \tx@NArray \psdots@iii}%
\end@SpecialObj}
\def\psdots@iii{%
\psk@dotsize
\@nameuse{psds@\psk@dotstyle}
newpath
n { transform floor .5 add exch floor .5 add exch itransform Dot } repeat}
\def\tx@SQ{SQ }
\def\tx@ST{ST }
\def\tx@SP{SP }
\def\pst@gdot#1{/Dot { gsave T \psk@dotangle \psk@dotscale #1 grestore } def }
\@namedef{psds@*}{\pst@gdot{0 0 DS \tx@SD}}
\@namedef{psds@o}{%
/r2 DS CLW sub def
\pst@gdot{0 0 DS \tx@SD \pst@usecolor\psfillcolor 0 0 r2 \tx@SD}}
\@namedef{psds@square*}{%
/r1 DS .886 mul def
\pst@gdot{r1 \tx@SQ}}
\@namedef{psds@square}{%
/r1 DS .886 mul def /r2 r1 CLW sub def
\pst@gdot{r1 \tx@SQ \pst@usecolor\psfillcolor r2 \tx@SQ}}
\@namedef{psds@triangle*}{%
/y1 DS .778 mul neg def /x1 y1 1.732 mul neg def
\pst@gdot{x1 y1 \tx@ST}}
\@namedef{psds@triangle}{%
/y1 DS .778 mul neg def /x1 y1 1.732 mul neg def
/y2 y1 CLW add def /x2 y2 1.732 mul neg def
\pst@gdot{x1 y1 \tx@ST \pst@usecolor\psfillcolor x2 y2 \tx@ST}}
\@namedef{psds@pentagon*}{%
/r1 DS 1.149 mul def
\pst@gdot{r1 \tx@SP}}
\@namedef{psds@pentagon}{%
DS .93 mul dup 1.236 mul /r1 ED CLW sub 1.236 mul /r2 ED
\pst@gdot{r1 \tx@SP \pst@usecolor\psfillcolor r2 \tx@SP}}
\@namedef{psds@+}{%
/DS DS 1.253 mul def
\pst@gdot{DS 0 moveto DS neg 0 L stroke 0 DS moveto 0 DS neg L stroke}}
\@namedef{psds@|}{%
\psk@tbarsize CLW mul add 2 div /DS ED
\pst@gdot{0 DS moveto 0 DS neg L stroke}}
\def\psset@dotstyle#1{%
\@ifundefined{psds@#1}%
{\@pstrickserr{Dot style `#1' not defined}\@eha}%
{\edef\psk@dotstyle{#1}}}
\psset@dotstyle{*}
\def\tx@FontDot{FontDot }
\def\newpsfontdot#1[#2]#3#4{%
  \@namedef{psds@#1}{%
    /#3 \psk@@dotangle [#2] \tx@FontDot
   /Dot { moveto gsave \psk@dotscale #4 show grestore } bind def 
  }%
}
\def\newpsfontdotH#1[#2]#3#4#5{%
  \@namedef{psds@#1}{%
    /#3 \psk@@dotangle [#2] \tx@FontDot
    /Dot {
      moveto
      \iftrue
      gsave \psk@dotscale \pst@usecolor\psfillcolor #5 show grestore
      \fi
      gsave \psk@dotscale #4 show grestore
    } bind def 
  }%
}
\pstheader{pst-dots.pro}
\newpsfontdot{*}[1.0 0.0 0.0 1.0 0.0 0.0]{PSTricksDotFont}{(b)}
\newpsfontdotH{o}[1.0 0.0 0.0 1.0 0.0 0.0]{PSTricksDotFont}{(c)}{(b)}
\newpsfontdotH{Bo}[1.0 0.0 0.0 1.0 0.0 0.0]{PSTricksDotFont}{(C)}{(b)}
\newpsfontdotH{triangle}[1.0 0.0 0.0 1.0 0.0 0.0]{PSTricksDotFont}{(t)}{(u)}
\newpsfontdotH{Btriangle}[1.0 0.0 0.0 1.0 0.0 0.0]{PSTricksDotFont}{(T)}{(u)}
\newpsfontdot{triangle*}[1.0 0.0 0.0 1.0 0.0 0.0]{PSTricksDotFont}{(u)}
\newpsfontdotH{square}[1.0 0.0 0.0 1.0 0.0 0.0]{PSTricksDotFont}{(s)}{(r)}
\newpsfontdotH{Bsquare}[1.0 0.0 0.0 1.0 0.0 0.0]{PSTricksDotFont}{(S)}{(r)}
\newpsfontdot{square*}[1.0 0.0 0.0 1.0 0.0 0.0]{PSTricksDotFont}{(r)}
\newpsfontdotH{pentagon}[1.0 0.0 0.0 1.0 0.0 0.0]{PSTricksDotFont}{(p)}{(q)}
\newpsfontdotH{Bpentagon}[1.0 0.0 0.0 1.0 0.0 0.0]{PSTricksDotFont}{(P)}{(q)}
\newpsfontdot{pentagon*}[1.0 0.0 0.0 1.0 0.0 0.0]{PSTricksDotFont}{(q)}
\newpsfontdotH{diamond}[1.0 0.0 0.0 1.0 0.0 0.0]{PSTricksDotFont}{(d)}{(l)}
\newpsfontdotH{Bdiamond}[1.0 0.0 0.0 1.0 0.0 0.0]{PSTricksDotFont}{(D)}{(l)}
\newpsfontdot{diamond*}[1.0 0.0 0.0 1.0 0.0 0.0]{PSTricksDotFont}{(l)}
\newpsfontdot{oplus}[1.44928 0.0 0.0 1.44928 -0.562319 -0.478261]{Symbol}{<C5>}
\newpsfontdot{otimes}[1.44928 0.0 0.0 1.44928 -0.562319 -0.475362]{Symbol}{<C4>}
\newpsfontdot{x}[1.8 0.0 0.0 1.8 -0.495 -0.4788]{Symbol}{<B4>}
\newpsfontdot{+}[2.3 0.0 0.0 2.3 -0.6486 -0.5819]{Times-Roman}{<2B>}
\newpsfontdot{asterisk}[2.43309 0.0 0.0 2.43309 -0.609489 -1.14477]{Times-Roman}{<2A>}
\newpsfontdot{B+}[2.3 0.0 0.0 2.3 -0.6555 -0.5819]{Times-Bold}{<2B>}
\newpsfontdot{Basterisk}[2.29358 0.0 0.0 2.29358 -0.576835 -1.08486]{Times-Bold}{<2A>}
\newpsfontdot{|}%
[1.98413 0.0 0.0 1.38 -0.258929 -0.5]{Helvetica}{(|)}
\newpsfontdot{B|}%
[1.98413 0.0 0.0 1.38 -0.277778 -0.5]{Helvetica-Bold}{(|)}
\newdimen\pslinearc
\def\psset@linearc#1{\pssetlength\pslinearc{#1}}
\psset@linearc{0pt}
\def\psline{\pst@object{psline}}
\def\psline@i{%
  \pst@getarrows{%
    \begin@OpenObj
    \pst@getcoors[\psline@ii%
  }%
}
\def\psline@ii{%
  \addto@pscode{\pst@cp \psline@iii \tx@Line}%
  \end@OpenObj%
}
\def\psline@iii{%
  \ifdim\pslinearc>\z@
    /r \pst@number\pslinearc def
    /Lineto { \tx@Arcto } def
  \else
    /Lineto /lineto load def
  \fi
  \ifshowpoints true \else false \fi
}
\def\qline(#1)(#2){%
  \def\pst@par{}%
  \begin@SpecialObj
  \def\pst@linetype{0}%
  \pst@getcoor{#1}\pst@tempa
  \pst@@getcoor{#2}%
  \addto@pscode{%
    \pst@tempa moveto \pst@coor L
    \@nameuse{psls@\pslinestyle}%
  }%
  \end@SpecialObj}
\def\pspolygon{\pst@object{pspolygon}}
\def\pspolygon@i{%
  \begin@ClosedObj%
  \def\pst@cp{}%
  \pst@getcoors[\pspolygon@ii%
}
\def\pspolygon@ii{%
  \addto@pscode{\psline@iii \tx@Polygon}%
  \def\pst@linetype{1}%
  \end@ClosedObj%
}
\def\psset@framearc#1{\pst@checknum{#1}\psk@framearc}
\psset@framearc{0}
\def\psset@cornersize#1{%
\pst@expandafter\psset@@cornersize{#1}\@nil}
\def\psset@@cornersize#1#2\@nil{%
\if #1a\relax
\def\psk@cornersize{\pst@number\pslinearc false }%
\else
\def\psk@cornersize{\psk@framearc true }%
\fi}
\psset@cornersize{relative}
\def\tx@Rect{Rect }
\def\tx@OvalFrame{OvalFrame }
\def\tx@Frame{Frame }
\def\psset@dimen#1{\pst@expandafter\psset@@dimen{#1}\@nil}
\def\psset@@dimen#1#2\@nil{%
  \if #1o\relax
    \def\psk@dimen{.5 }%
  \else
    \if #1m\relax
      \def\psk@dimen{0 }%
    \else
      \if #1i\relax
        \def\psk@dimen{-.5 }%
  \fi\fi\fi}
\psset@dimen{outer}
\def\psframe{\pst@object{psframe}}
\def\psframe@i(#1){%
  \@ifnextchar({\psframe@ii(#1)}{\psframe@ii(0,0)(#1)}}
\def\psframe@ii(#1)(#2){%
  \begin@ClosedObj
    \pst@getcoor{#1}\pst@tempa
    \pst@@getcoor{#2}%
    \addto@pscode{\psk@cornersize \pst@tempa \pst@coor \psk@dimen \tx@Frame}%
    \def\pst@linetype{2}%
    \showpointsfalse
  \end@ClosedObj%
}
\def\tx@BezierNArray{BezierNArray }
\def\tx@OpenBezier{OpenBezier }
\def\tx@ClosedBezier{ClosedBezier }
\def\tx@BezierShowPoints{BezierShowPoints }
\def\psbezier{\pst@object{psbezier}}
\def\psbezier@i{%
  \pst@getarrows{%
    \begin@OpenObj
      \pst@getcoors[\psbezier@ii%
  }%
}
\def\psbezier@ii{%
  \addto@pscode{%
    \ifshowpoints true \else false \fi
    \tx@OpenBezier
    \ifshowpoints \tx@BezierShowPoints \fi}%
  \end@OpenObj}
\def\pscbezier{\pst@object{pscbezier}}
\def\pscbezier@i{%
  \begin@ClosedObj
  \pst@getcoors[\pscbezier@ii}
\def\pscbezier@ii{%
  \addto@pscode{%
    \ifshowpoints true \else false \fi
    \tx@ClosedBezier
    \ifshowpoints \tx@BezierShowPoints \fi}%
  \chardef\pst@linetype=1
  \end@ClosedObj}
\def\tx@Parab{Parab }
\def\parabola{\pst@object{parabola}}
\def\parabola@i{\pst@getarrows\parabola@ii}
\def\parabola@ii#1(#2)#3(#4){%
\begin@OpenObj
\pst@getcoor{#2}\pst@tempa
\pst@@getcoor{#4}%
\addto@pscode{\pst@tempa \pst@coor \tx@Parab}%
\end@OpenObj}
\def\psset@gridwidth#1{\pst@getlength{#1}\psk@gridwidth}
\psset@gridwidth{.8pt}
\def\psset@griddots#1{%
  \pst@cntg=#1\relax
  \edef\psk@griddots{\the\pst@cntg}}
\psset@griddots{0}
\def\psset@gridcolor#1{\pst@getcolor{#1}\psgridcolor}
\psset@gridcolor{black}
\def\psset@subgridwidth#1{\pst@getlength{#1}\psk@subgridwidth}
\psset@subgridwidth{.4pt}
\def\psset@subgridcolor#1{\pst@getcolor{#1}\pssubgridcolor}
\psset@subgridcolor{gray}
\def\psset@subgriddots#1{%
  \pst@cntg=#1\relax\edef\psk@subgriddots{\the\pst@cntg}}
\psset@subgriddots{0}
\def\psset@subgriddiv#1{%
  \pst@cntg=#1\relax\edef\psk@subgriddiv{\the\pst@cntg}}
\psset@subgriddiv{5}
\def\psset@gridlabels#1{\pst@getlength{#1}\psk@gridlabels}
\psset@gridlabels{10pt}
\def\psset@gridlabelcolor#1{\pst@getcolor{#1}\psgridlabelcolor}
\psset@gridlabelcolor{black}
\def\tx@Grid{Grid }
\def\psgrid{\pst@object{psgrid}}
\def\psgrid@i{\@ifnextchar(%
  {\psgrid@ii}{\expandafter\psgrid@iv\pic@coor}}
\def\psgrid@ii(#1){\@ifnextchar(%
  {\psgrid@iii(#1)}{\psgrid@iv(0,0)(0,0)(#1)}}
\def\psgrid@iii(#1)(#2){\@ifnextchar(%
  {\psgrid@iv(#1)(#2)}{\psgrid@iv(#1)(#1)(#2)}}
\def\psgrid@iv(#1)(#2)(#3){%
  \begin@SpecialObj
    \pst@getcoor{#1}\pst@tempA
    \pst@getcoor{#2}\pst@tempB 
    \pst@@getcoor{#3}%
    \ifnum\psk@subgriddiv>1
      \addto@pscode{
        gsave
        \psk@subgridwidth SLW \pst@usecolor\pssubgridcolor
        \pst@tempB \pst@coor \pst@tempA                 
        \pst@number\psxunit abs \pst@number\psyunit abs 
        \psk@subgriddiv\space \psk@subgriddots\space
        {} 0 \tx@Grid 
	grestore
      }%
    \fi
    \addto@pscode{%
      gsave
      \psk@gridwidth SLW \pst@usecolor\psgridcolor
      \pst@tempB \pst@coor \pst@tempA                 
      \pst@number\psxunit abs \pst@number\psyunit abs 
      1 \psk@griddots\space { \pst@usecolor\psgridlabelcolor }
      \psk@gridlabels \tx@Grid 
      grestore
    }%
  \end@SpecialObj%
}
\newif\ifpsmathbox
\psmathboxtrue
\def\pst@mathflag{\z@}
\newtoks\everypsbox
\let\pst@thisbox\relax
\long\def\pst@makenotverbbox#1#2{%
  \edef\pst@mathflag{%
  \ifpsmathbox\ifmmode\ifinner 1\else 2\fi\else \z@\fi\else \z@\fi}%
  \setbox\pst@hbox=\hbox{%
    \ifcase\pst@mathflag\or$\m@th\textstyle\or$\m@th\displaystyle\fi
    {\pst@thisbox\the\everypsbox#2}%
    \ifnum\pst@mathflag>\z@$\fi
  }%
  #1%
}
\def\pst@makeverbbox#1{%
  \def\pst@afterbox{#1}%
  \edef\pst@mathflag{%
    \ifpsmathbox\ifmmode\ifinner 1\else 2\fi\else \z@\fi\else \z@\fi%
  }%
  \afterassignment\pst@beginbox
  \setbox\pst@hbox\hbox%
}
\def\pst@beginbox{%
  \ifcase\pst@mathflag\or$\m@th\or$\m@th\displaystyle\fi
  \bgroup\aftergroup\pst@endbox
  \pst@thisbox
  \the\everypsbox%
}
\def\pst@endbox{%
  \ifnum\pst@mathflag>\z@$\fi
  \egroup
  \pst@afterbox%
}
\def\pst@makebox{\pst@@makebox}
\def\psverbboxtrue{\def\pst@@makebox{\pst@makeverbbox}}
\def\psverbboxfalse{\def\pst@@makebox{\pst@makenotverbbox}}
\psverbboxfalse
\def\pst@longbox{%
  \def\pst@makebox{%
    \gdef\pst@makebox{\pst@@makebox}%
    \pst@makelongbox%
  }%
}
\def\pst@makelongbox#1{%
  \def\pst@afterbox{#1}%
  \edef\pst@mathflag{%
    \ifpsmathbox\ifmmode\ifinner 1\else 2\fi\else \z@\fi\else \z@\fi%
  }%
  \setbox\pst@hbox\hbox\bgroup
  \aftergroup\pst@afterbox
  \ifcase\pst@mathflag\or$\m@th\or$\m@th\displaystyle\fi
  \begingroup
  \pst@thisbox
  \the\everypsbox%
}
\def\pst@endlongbox{%
  \endgroup
  \ifnum\pst@mathflag>\z@$\fi
  \egroup%
}
\def\pslongbox#1#2{%
\@namedef{#1}{\pst@longbox#2}%
\@namedef{end#1}{\pst@endlongbox}%
}
\newdimen\psframesep
\def\psset@framesep#1{\pssetlength\psframesep{#1}}
\psset@framesep{3pt}
\newif\ifpsboxsep
\def\psset@boxsep#1{\@nameuse{psboxsep#1}}
\psset@boxsep{true}
\def\pst@useboxpar{%
  \use@par
  \if@star
    \let\pslinecolor\psfillcolor
    \solid@star
    \let\solid@star\relax
  \fi
  \ifpsdoubleline \pst@setdoublesep \fi%
}
\def\psframebox{\pst@object{psframebox}}
\def\psframebox@i{\pst@makebox\psframebox@ii}
\def\psframebox@ii{%
  \begingroup
  \pst@useboxpar
  \pst@dima=\pslinewidth
  \advance\pst@dima by \psframesep
  \pst@dimc=\wd\pst@hbox\advance\pst@dimc by \pst@dima
  \pst@dimb=\dp\pst@hbox\advance\pst@dimb by \pst@dima
  \pst@dimd=\ht\pst@hbox\advance\pst@dimd by \pst@dima
  \setbox\pst@hbox=\hbox{%
    \ifpsboxsep\kern\pst@dima\fi
    \begin@ClosedObj
    \addto@pscode{%
      \psk@cornersize
      \pst@number\pst@dima neg
      \pst@number\pst@dimb neg
      \pst@number\pst@dimc
      \pst@number\pst@dimd
      .5
      \tx@Frame%
    }%
    \def\pst@linetype{2}%
    \showpointsfalse
    \end@ClosedObj
    \box\pst@hbox
    \ifpsboxsep\kern\pst@dima\fi%
  }%
  \ifpsboxsep\dp\pst@hbox=\pst@dimb\ht\pst@hbox=\pst@dimd\fi
  \leavevmode\box\pst@hbox
  \endgroup%
}
\def\psdblframebox{\pst@object{psdblframebox}}
\def\psdblframebox@i{\addto@par{doubleline=true}\psframebox@i}
\def\psclip#1{%
  \leavevmode
  \begingroup
    \begin@psclip
      \begingroup
        \def\use@pscode{%
          \pstVerb{%
            \pst@dict
            /mtrxc CM def
            CP CP T
            \tx@STV
            \psk@origin
            \psk@swapaxes
            newpath
            \pst@code
            clip
            newpath
            mtrxc setmatrix
            moveto
            0 setgray
            end
	  }%
          \gdef\pst@code{}}%
  \def\@multips(##1)(##2)##3##4{\pst@misplaced\multips}%
  \def\nc@object##1##2##3##4{\pst@misplaced{node connection}}%
  \hbox to\z@{#1}%
  \endgroup
  \def\endpsclip{%
    \end@psclip
    \endgroup}%
  \ignorespaces}
\def\endpsclip{\pst@misplaced\endpsclip}
\let\begin@psclip\relax
\def\end@psclip{\pstVerb{currentpoint initclip moveto}}
\def\AltClipMode{%
\def\end@psclip{\pstVerb{\pst@grestore}}%
\def\begin@psclip{\pstVerb{gsave}}}
\def\clipbox{\@ifnextchar[{\clipbox@}{\clipbox@[\z@]}}
\def\clipbox@[#1]{\pst@makebox{\clipbox@@{#1}}}
\def\clipbox@@#1{%
  \pssetlength\pst@dimg{#1}%
  \leavevmode\hbox{%
  \begin@psclip
  \pst@Verb{%
    CM \tx@STV CP T newpath
    /a \pst@number\pst@dimg def
    /w \pst@number{\wd\pst@hbox}a add def
    /d \pst@number{\dp\pst@hbox}a add neg def
    /h \pst@number{\ht\pst@hbox}a add def
    a neg d moveto
    a neg h L
    w h L
    w d L
    closepath
    clip
    newpath
    0 0 moveto
    setmatrix}%
\unhbox\pst@hbox
\end@psclip}}
\def\psshadowbox{\pst@object{psshadowbox}}
\def\psshadowbox@i{\pst@makebox\psshadowbox@ii}
\def\psshadowbox@ii{%
\begingroup
\pst@useboxpar
\psshadowtrue
\psboxseptrue
\def\psk@shadowangle{-45 }%
\setbox\pst@hbox=\hbox{\psframebox@ii}%
\pst@dimh=\psk@shadowsize\p@
\pst@dimh=.7071\pst@dimh
\pst@dimg=\dp\pst@hbox
\advance\pst@dimg\pst@dimh
\dp\pst@hbox=\pst@dimg
\pst@dimg=\wd\pst@hbox
\advance\pst@dimg\pst@dimh
\wd\pst@hbox=\pst@dimg
\leavevmode
\box\pst@hbox
\endgroup}
\def\pscirclebox{\pst@object{pscirclebox}}
\def\pscirclebox@i{\pst@makebox\pscirclebox@ii}
\def\pscirclebox@ii{%
\begingroup
\pst@useboxpar
\setbox\pst@hbox=\hbox{%
\pst@nodehook
\pscirclebox@iii
\box\pst@hbox}%
\ifpsboxsep \pscirclebox@sep \fi
\leavevmode
\box\pst@hbox
\endgroup}
\def\pscirclebox@iii{%
\if@star
\pslinewidth\z@
\pstverb{\pst@dict \tx@STP \pst@usecolor\psfillcolor
newpath \pscirclebox@iv \tx@SD end}%
\else
\begin@ClosedObj
\def\pst@linetype{4}\showpointsfalse
\addto@pscode{%
\pscirclebox@iv CLW 2 div add 0 360 arc closepath}%
\end@ClosedObj
\fi}
\def\pscirclebox@iv{%
\pst@number{\wd\pst@hbox}2 div
\pst@number{\ht\pst@hbox}\pst@number{\dp\pst@hbox}add 2 div
2 copy \pst@number{\dp\pst@hbox}sub 4 2 roll
\tx@Pyth \pst@number\psframesep add }
\def\pscirclebox@sep{%
\pst@dimb=\ht\pst@hbox
\advance\pst@dimb\dp\pst@hbox
\divide\pst@dimb 2
\pst@dima=.5\wd\pst@hbox
\pst@pyth\pst@dima\pst@dimb\pst@dimc
\advance\pst@dimc\pslinewidth
\advance\pst@dimc\psframesep
\advance\pst@dimb-\pst@dimc
\setbox\pst@hbox=\hbox to2\pst@dimc{%
\hss
\vbox{\kern-\pst@dimb\box\pst@hbox}%
\hss}%
\advance\pst@dimb-\dp\pst@hbox
\dp\pst@hbox=-\pst@dimb}
\let\pst@nodehook\relax
\def\psCirclebox{\pst@object{psCirclebox}}
\def\psCirclebox@i{\pst@makebox\psCirclebox@ii}
\def\psCirclebox@ii{%
\begingroup
\pst@useboxpar
\pst@dima=\ht\pst@hbox
\advance\pst@dima\dp\pst@hbox
\divide\pst@dima\tw@
\pssetlength\pst@dimb\psk@radius
\setbox\pst@hbox=\hbox{%
\pst@nodehook
\pscircle(.5\wd\pst@hbox,\pst@dima){\pst@dimb}%
\box\pst@hbox}%
\ifpsboxsep \psCirclebox@sep \fi
\leavevmode
\box\pst@hbox
\endgroup}
\def\psCirclebox@sep{%
\pst@dimc=\pst@dimb
\advance\pst@dimb-\pst@dima
\advance\pst@dima\pst@dimc
\setbox\pst@hbox=\hbox to\tw@\pst@dimc{%
\hss
\vrule width \z@ depth \pst@dimb height \pst@dima
\box\pst@hbox
\hss}}%
\def\psovalbox{\pst@object{psovalbox}}
\def\psovalbox@i{\pst@makebox{\psovalbox@ii}}
\def\psovalbox@ii{%
\begingroup
\pst@useboxpar
\psovalbox@iii
\ifpsboxsep\psovalbox@sep\fi
\leavevmode
\box\pst@hbox
\endgroup}
\def\psovalbox@iii{%
\psovalbox@iv
\setbox\pst@hbox=\hbox{%
\begin@ClosedObj
\addto@pscode{%
0 360
\pst@number\pst@dimc CLW 2 div sub
\pst@number\pst@dimd CLW 2 div sub
\pst@number\pst@dima
\pst@number\pst@dimb
\tx@Ellipse
closepath}%
\def\pst@linetype{2}%
\end@ClosedObj
\unhbox\pst@hbox}}
\def\psovalbox@iv{%
\pst@dimc=\pslinewidth\advance\pst@dimc\psframesep
\pst@dimd=\ht\pst@hbox\advance\pst@dimd\dp\pst@hbox
\pst@dima=.5\wd\pst@hbox
\pst@dimb=.5\pst@dimd\advance\pst@dimb-\dp\pst@hbox
\pst@dimd=.707\pst@dimd
\advance\pst@dimd\pst@dimc
\advance\pst@dimc.707\wd\pst@hbox}
\def\psovalbox@sep{%
\setbox\pst@hbox\hbox to 2\pst@dimc{\hss\unhbox\pst@hbox\hss}%
\pst@dimg=\pst@dimd
\advance\pst@dimg-\pst@dimb
\dp\pst@hbox=\pst@dimg
\advance\pst@dimd\pst@dimb
\ht\pst@hbox=\pst@dimd}
\def\psdiabox{\pst@object{psdiabox}}
\def\psdiabox@i{\pst@makebox{\psdiabox@ii}}
\def\psdiabox@ii{%
\begingroup
\pst@useboxpar
\psdiabox@iii
\ifpsboxsep\psdiabox@sep\fi
\leavevmode
\box\pst@hbox
\endgroup}
\def\psdiabox@iv{%
\pst@dimg=.707\pslinewidth
\advance\pst@dimg.707\psframesep
\pst@dima=\wd\pst@hbox
\divide\pst@dima 2
\pst@dimc=\pst@dima
\advance\pst@dimc\pst@dimg
\pst@dimd=\ht\pst@hbox
\advance\pst@dimd\dp\pst@hbox
\divide\pst@dimd 2
\pst@dimb=\pst@dimd
\advance\pst@dimb-\dp\pst@hbox
\advance\pst@dimd\pst@dimg}
\def\psdiabox@iii{%
\psdiabox@iv
\setbox\pst@hbox=\hbox{%
\begin@ClosedObj
\addto@pscode{%
\psline@iii
pop
.5
\pst@number\pst@dimc 2 mul \pst@number\pst@dimd 2 mul
0
\pst@number\pst@dima \pst@number\pst@dimb
\tx@Diamond}%
\def\pst@linetype{4}%
\end@ClosedObj
\box\pst@hbox}}
\def\psdiabox@sep{%
\setbox\pst@hbox\hbox to 4\pst@dimc{\hss\unhbox\pst@hbox\hss}%
\multiply\pst@dimd 2
\advance\pst@dimd\pst@dimb
\ht\pst@hbox\pst@dimd
\advance\pst@dimd-2\pst@dimb
\dp\pst@hbox\pst@dimd}
\def\psset@trimode#1{\pst@expandafter\psset@@trimode{#1}\@empty\@empty\@nil}
\def\psset@@trimode#1#2#3\@nil{%
\let\pst@tempg#1\relax
\ifx\pst@tempg*%
\let\psk@@trimode\@empty
\let\pst@tempg#2\relax
\else
\let\psk@@trimode\relax
\fi
\edef\psk@trimode{%
\ifx R\pst@tempg 1 \else\ifx D\pst@tempg 2
\else\ifx L\pst@tempg 3 \else 0 \fi\fi\fi}}
\psset@trimode{U}
\def\pstribox{\pst@object{pstribox}}
\def\pstribox@i{\pst@makebox{\pstribox@ii}}
\def\pstribox@ii{%
\begingroup
\pst@useboxpar
\pstribox@iii
\ifpsboxsep\pstribox@sep\fi
\leavevmode
\box\pst@hbox
\endgroup}
\def\pstribox@iii{%
\pstribox@iv
\setbox\pst@hbox=\hbox{%
\begin@ClosedObj
\addto@pscode{%
\psline@iii
pop
.5
\pst@number\pst@dimc \pst@number\pst@dimd
\ifodd\psk@trimode exch \fi
\psk@trimode -90 mul
\pst@number\pst@dima \pst@number\pst@dimb
\tx@Triangle}%
\def\pst@linetype{2}%
\end@ClosedObj
\box\pst@hbox}}
\def\pstribox@iv{%
\pst@dimh=\pslinewidth
\advance\pst@dimh\psframesep
\pst@dimg=\ht\pst@hbox
\advance\pst@dimg-\dp\pst@hbox
\divide\pst@dimg 2
\edef\pst@tempa{\number\pst@dimg sp}
\ifodd\psk@trimode
\pst@dimb\pst@dimg
\else
\pst@dima=\wd\pst@hbox
\divide\pst@dima 2
\fi
\ifcase\psk@trimode
\pst@dimb=-\dp\pst@hbox
\advance\pst@dimb-\pst@dimh
\or
\pst@dima=-\pst@dimh
\or
\pst@dimb=\ht\pst@hbox
\advance\pst@dimb\pst@dimh
\or
\pst@dima=\wd\pst@hbox
\advance\pst@dima\pst@dimh
\fi
\pst@dimd=\dp\pst@hbox
\advance\pst@dimd\ht\pst@hbox
\ifx\psk@@trimode\relax
\pst@dimc=\wd\pst@hbox
\advance\pst@dimc\ifodd\psk@trimode 1.447\else 1.789\fi\pst@dimh
\multiply\pst@dimc 2
\advance\pst@dimd\ifodd\psk@trimode 1.789\else 1.447\fi\pst@dimh
\multiply\pst@dimd 2
\else
\ifodd\psk@trimode
\advance\pst@dimd 1.1547\wd\pst@hbox
\advance\pst@dimd 3.4641\pst@dimh
\pst@dimc=.866\pst@dimd
\else
\advance\pst@dimd .866\wd\pst@hbox 
\advance\pst@dimd 3\pst@dimh
\pst@dimc=1.1547\pst@dimd 
\fi
\fi}
\def\pstribox@sep{%
\ifodd\psk@trimode
\advance\pst@dimb.5\pst@dimd
\ht\pst@hbox=\pst@dimb
\advance\pst@dimd-\pst@dimb
\dp\pst@hbox=\pst@dimd
\else
\setbox\pst@hbox\hbox to \pst@dimc{\hss\unhbox\pst@hbox\hss}%
\global\pst@dimg=.5\pst@dimc
\fi
\ifcase\psk@trimode
\dp\pst@hbox-\pst@dimb
\advance\pst@dimd\pst@dimb
\ht\pst@hbox\pst@dimd
\or
\pst@dimg=.5\wd\pst@hbox
\global\advance\pst@dimg-\pst@dima
\setbox\pst@hbox\hbox to \pst@dimc{\kern-\pst@dima\box\pst@hbox\hss}%
\or
\ht\pst@hbox\pst@dimb
\advance\pst@dimd-\pst@dimb
\dp\pst@hbox\pst@dimd
\or
\pst@dimg=\pst@dimc
\advance\pst@dimg-\pst@dima
\global\advance\pst@dimg.5\wd\pst@hbox
\setbox\pst@hbox\hbox to \pst@dimc{%
\hss\box\pst@hbox\kern\psframesep\kern\pslinewidth}%
\fi}
\def\psset@arcsepA#1{\pst@getlength{#1}\psk@arcsepA}
\def\psset@arcsepB#1{\pst@getlength{#1}\psk@arcsepB}
\def\psset@arcsep#1{%
\psset@arcsepA{#1}\let\psk@arcsepB\psk@arcsepA}
\psset@arcsep{0}
\def\tx@ArcArrow{ArcArrow }
\def\psarc{\pst@object{psarc}}
\def\psarc@i{\@ifnextchar({\psarc@iii}{\psarc@ii}}
\def\psarc@ii#1{\addto@par{arrows=#1}%
  \@ifnextchar({\psarc@iii}{\psarc@iii(0,0)}%
}
\def\psarc@iii(#1)#2#3#4{%
  \begin@OpenObj
    \pst@getangle{#3}\pst@tempa
    \pst@getangle{#4}\pst@tempb
    \pst@@getcoor{#1}%
    \pssetlength\pst@dima{#2}%
    \addto@pscode{\psarc@iv \psarc@v}%
    \gdef\psarc@type{0}%
    \showpointsfalse
  \end@OpenObj%
}
\def\psarc@iv{%
  \pst@coor /y ED /x ED
  /r \pst@number\pst@dima def
  /c 57.2957 r \tx@Div def
  /angleA
    \pst@tempa
    \psk@arcsepA c mul 2 div
    \ifcase \psarc@type add \or sub \fi
  def
  /angleB
    \pst@tempb
    \psk@arcsepB c mul 2 div
    \ifcase \psarc@type sub \or add \fi
  def
  \ifshowpoints\psarc@showpoints\fi
  \ifx\psk@arrowA\@empty
    \ifnum\psk@liftpen=2
      r angleA \tx@PtoC
      y add exch x add exch moveto
    \fi
  \fi}
\def\psarc@v{%
  x y r
  angleA
  \ifx\psk@arrowA\@empty\else
    { ArrowA CP }
    r 0 gt \pslbrace
    { \ifcase\psarc@type add \or sub \fi } \psrbrace\pslbrace
    { \ifcase\psarc@type sub \or add \fi } \psrbrace ifelse
    \tx@ArcArrow
  \fi
  angleB
  \ifx\psk@arrowB\@empty\else
    { ArrowB }
    r 0 gt \pslbrace
      { \ifcase\psarc@type sub \or add \fi } \psrbrace\pslbrace
      { \ifcase\psarc@type add \or sub \fi } \psrbrace ifelse
 \tx@ArcArrow
  \fi
\ifcase\psarc@type arc \or arcn \fi}
%
\def\psarc@type{0}
\def\psarc@showpoints{%
  gsave
  newpath
  x y moveto
  x y r \pst@tempa \pst@tempb
  \ifcase\psarc@type arc \or arcn \fi
  closepath
  CLW 2 div SLW
  [ \psk@dash\space ] 0 setdash stroke
  grestore }
\def\psarcn{\pst@object{psarcn}}
\def\psarcn@i{\def\psarc@type{1}\psarc@i}
%
%
%
\def\psellipticwedge{\def\pst@par{}\pst@object{psellipticwedge}}
\def\psellipticwedge@i(#1){%
  \@ifnextchar({\psellipticwedge@ii(#1)}{\psellipticwedge@ii(0,0)(#1)}}
\def\psellipticwedge@ii(#1)(#2)#3#4{%
  \begin@ClosedObj
    \pst@getangle{#3}\pst@tempa
    \pst@getangle{#4}\pst@tempb
    \pst@getcoor{#1}\pst@tempc
    \pst@@getcoor{#2}%
    \def\pst@linetype{1}%
    \addto@pscode{%
      \pst@tempa \pst@tempb
      \pst@coor
      \pst@tempc moveto
      \ifdim\psk@dimen\p@=\z@\else
        \psk@dimen CLW mul dup 3 1 roll
        sub 3 1 roll sub exch
      \fi
      \pst@tempc
      \tx@Ellipse
      closepath%
    }%
    \showpointsfalse
  \end@ClosedObj%
}
%
%
\def\psellipticarcn{\def\pst@par{}\pst@object{psellipticarcn}}
\def\psellipticarcn@i{\let\if@psarcn\iftrue\psellipticarc@ii}
\def\psellipticarc{\def\pst@par{}\pst@object{psellipticarc}}
\def\psellipticarc@i{\let\if@psarcn\iffalse\psellipticarc@ii}
\let\if@psarcn\iffalse
\def\psellipticarc@ii{\pst@getarrows\psellipticarc@iii}
\def\psellipticarc@iii(#1){%
	\@ifnextchar({\psellipticarc@iv(#1)}{\psellipticarc@iv(0,0)(#1)}}
\def\psellipticarc@iv(#1)(#2)#3#4{%
  \begin@OpenObj
  \pst@getcoor{#1}\pst@tempa
  \pst@getcoor{#2}\pst@tempb
  \pst@getangle{#3}\pst@tempc
  \pst@getangle{#4}\pst@tempd
  \addto@pscode{\psellipticarc@definearg \psellipticarc@draw}%
  \showpointsfalse
  \end@OpenObj%
}
\def\psellipticarc@definearg{%
  \pst@tempa /y ED /x ED  
  \pst@tempb              
  \ifdim\psk@dimen\p@=\z@\else
    \psk@dimen CLW mul dup 3 1 roll
    sub 3 1 roll sub exch
  \fi
  /ry ED /rx ED
  /angleA
    /d {  \if@psarcn sub \else add \fi } def
    \pst@tempc \psk@arcsepA 2 div
    ArcAdjust
  def
  /angleB
    /d  {  \if@psarcn add \else sub \fi } def
      \pst@tempd \psk@arcsepB 2 div
      ArcAdjust
    def
    \ifshowpoints\psellipticarc@showpoints\fi
    \ifx\psk@arrowA\@empty
      \ifnum\psk@liftpen=2
	angleA cos rx mul x add
	angleA sin ry mul y add
	moveto
     \fi
   \fi%
}
\def\psellipticarc@draw{%
  0 0 1
  angleA
  \ifx\psk@arrowA\@empty\else
    { ArrowA CP }
    { \if@psarcn sub \else add \fi }
    EllipticArcArrow
  \fi
  angleB
  \ifx\psk@arrowB\@empty\else
    { ArrowB }
    { \if@psarcn add \else sub \fi }
    EllipticArcArrow
  \fi
  /mtrx CM def
  x y T
  rx ry scale
  \if@psarcn arcn \else arc \fi
  mtrx setmatrix%
}
\def\psellipticarc@showpoints{%
  gsave
  /mtrx CM def
  x y T
  rx ry scale
  0 0 moveto
  0 0 1 \pst@tempc \pst@tempd
  \ifcase\psarc@type arc \or arcn \fi
  closepath
  mtrx setmatrix
  CLW 2 div SLW
  [ \psk@dash\space ] 0 setdash stroke
  grestore %
}
\def\pscircle{\pst@object{pscircle}}
\def\pscircle@i{\@ifnextchar({\pscircle@do}{\pscircle@do(0,0)}}
\def\pscircle@do(#1)#2{%
\if@star
{\use@par\qdisk(#1){#2}}%
\else
\begin@ClosedObj
\pst@@getcoor{#1}%
\pssetlength\pst@dimc{#2}%
\def\pst@linetype{4}%
\addto@pscode{%
\pst@coor
\pst@number\pst@dimc
\psk@dimen CLW mul sub
0 360 arc
closepath}%
\showpointsfalse
\end@ClosedObj
\fi
\ignorespaces}
\def\qdisk(#1)#2{%
\def\pst@par{}%
\begin@SpecialObj
\pst@@getcoor{#1}%
\pssetlength\pst@dimg{#2}%
\addto@pscode{\pst@coor \pst@number\pst@dimg \tx@SD}%
\end@SpecialObj}
\def\psset@radius#1{\pst@@getlength{#1}\psk@radius}
\psset@radius{.25cm}
\def\psCircle{\pst@object{psCircle}}
\def\psCircle@i{\@ifnextchar({\psCircle@ii}{\psCircle@ii(0,0)}}
\def\psCircle@ii(#1){\pscircle@do(#1){\psk@radius}}
\def\pswedge{\pst@object{pswedge}}
\def\pswedge@i{\@ifnextchar({\pswedge@ii}{\pswedge@ii(0,0)}}
\def\pswedge@ii(#1)#2#3#4{%
  \begin@ClosedObj
  \pssetlength\pst@dimc{#2}
  \pst@getangle{#3}\pst@tempa
  \pst@getangle{#4}\pst@tempb
  \pst@@getcoor{#1}%
  \def\pst@linetype{1}%
  \addto@pscode{%
    \pst@coor
    2 copy
    moveto
    \pst@number\pst@dimc \psk@dimen CLW mul sub 
    \pst@tempa \pst@tempb
    arc
    closepath}%
  \showpointsfalse
  \end@ClosedObj%
}
\def\tx@Ellipse{Ellipse }
\def\psellipse{\pst@object{psellipse}}
\def\psellipse@i(#1){\@ifnextchar({\psellipse@ii(#1)}{\psellipse@ii(0,0)(#1)}}
\def\psellipse@ii(#1)(#2){%
  \begin@ClosedObj
  \pst@getcoor{#1}\pst@tempa
  \pst@@getcoor{#2}%
  \addto@pscode{%
    0 360
    \pst@coor
    \ifdim\psk@dimen\p@=\z@\else
      \psk@dimen CLW mul
      dup 4 -1 roll sub neg 3 1 roll sub
    \fi
    \pst@tempa
    \tx@Ellipse
    closepath%
  }%
  \def\pst@linetype{2}%
  \end@ClosedObj%
}
\def\multips{\@ifnextchar({\def\pst@par{}\multips@ii}{\multips@i}}
\def\multips@i#1{\def\pst@par{rot=#1}\multips@ii}
\def\multips@ii(#1){\@ifnextchar({\multips@iii(#1)}{\multips@iii(\z@,\z@)(#1)}}
\long\def\multips@iii(#1)(#2)#3#4{%
  \begingroup
  \pst@killglue
  \use@par
  \pst@getcoor{#1}\pst@tempa
  \pst@@getcoor{#2}%
  \pst@cnta=#3\relax
  \init@pscode
  \addto@pscode{%
    \pst@tempa T \the\pst@cnta\space \pslbrace
    gsave \ifx\psk@rot\@empty\else\psk@rot rotate \fi}%
  \hbox to\z@{%
    \def\init@pscode{%
      \addto@pscode{%
        gsave
        \pst@number\pslinewidth SLW
        \pst@usecolor\pslinecolor}}%
    \def\use@pscode{\addto@pscode{grestore}}%
    \def\psclip##1{\pst@misplaced\psclip}%
    \def\nc@object##1##2##3##4{\pst@misplaced{node connection}}%
    #4%
  }%
  \addto@pscode{grestore \pst@coor T \psrbrace repeat}%
  \leavevmode
  \use@pscode
  \endgroup
  \ignorespaces}
\def\psscalebox#1{\pst@makebox{\ps@scalebox{#1}}}
\def\ps@scalebox#1{%
  \begingroup
  \pst@getscale{#1}\pst@tempa
  \let\pst@tempc\pst@tempg
  \let\pst@tempd\pst@temph
  \ps@@scalebox
  \endgroup}
\def\ps@@scalebox{%
  \leavevmode
  \hbox{%
    \ifdim\pst@tempd\p@<\z@
      \pst@dimg=\pst@tempd\ht\pst@hbox
      \pst@dimh=\pst@tempd\dp\pst@hbox
      \dp\pst@hbox=-\pst@dimg
      \ht\pst@hbox=-\pst@dimh
    \else
      \ht\pst@hbox=\pst@tempd\ht\pst@hbox
      \dp\pst@hbox=\pst@tempd\dp\pst@hbox
    \fi
    \pst@dima=\pst@tempc\wd\pst@hbox
    \ifdim\pst@dima<\z@\kern-\pst@dima\fi
    \pst@Verb{CP CP translate \pst@tempa \tx@NET}%
    \hbox to \z@{\box\pst@hbox\hss}%
    \pst@Verb{%
      CP CP translate
      1 \pst@tempc div 1 \pst@tempd div scale
      \tx@NET}%
    \ifdim\pst@dima>\z@\kern\pst@dima\fi%
  }%
}
\pslongbox{Scalebox}{\psscalebox}
\def\psscaleboxto(#1,#2){\pst@makebox{\ps@scaleboxto(#1,#2)}}
\def\ps@scaleboxto(#1,#2){%
  \begingroup
  \pssetlength\pst@dima{#1}%
  \pssetlength\pst@dimb{#2}%
  \ifdim\pst@dima=\z@\else
    \pst@divide{\pst@dima}{\wd\pst@hbox}\pst@tempc
    \edef\pst@tempc{\pst@tempc\space}%
  \fi
  \ifdim\pst@dimb=\z@
    \ifdim\pst@dima=\z@
      \@pstrickserr{%
        \string\psscaleboxto\space dimensions cannot both be zero}\@ehpa
      \def\pst@tempa{}%
      \def\pst@tempc{1 }%
      \def\pst@tempd{1 }%
    \else
      \let\pst@tempd\pst@tempc
    \fi
  \else
    \pst@dimc=\ht\pst@hbox
    \advance\pst@dimc\dp\pst@hbox
    \pst@divide{\pst@dimb}{\pst@dimc}\pst@tempd
    \edef\pst@tempd{\pst@tempd\space}%
    \ifdim\pst@dima=\z@ \let\pst@tempc\pst@tempd \fi
  \fi
  \edef\pst@tempa{\pst@tempc \pst@tempd scale }%
  \ps@@scalebox
  \endgroup}
\pslongbox{Scaleboxto}{\psscaleboxto}
\def\tx@Rot{Rot }
\def\psrotateleft{\pst@makebox{\ps@rotateleft\pst@hbox}}
\def\ps@rotateleft#1{%
\leavevmode\hbox{\hskip\ht#1\hskip\dp#1\vbox{\vskip\wd#1%
\pst@Verb{90 \tx@Rot}
\vbox to \z@{\vss\hbox to \z@{\box#1\hss}\vskip\z@}%
\pst@Verb{-90 \tx@Rot}}}}
\def\psrotateright{\pst@makebox{\ps@rotateright\pst@hbox}}
\def\ps@rotateright#1{%
  \leavevmode\hbox{%
  \hskip\ht#1\hskip\dp#1\vbox{\vskip\wd#1%
  \pst@Verb{-90 \tx@Rot}
  \vbox to \z@{\hbox to \z@{\hss\box#1}\vss}%
  \pst@Verb{90 \tx@Rot}}}}
\def\psrotatedown{\pst@makebox{\ps@rotatedown\pst@hbox}}
\def\ps@rotatedown#1{%
\hbox{\hskip\wd#1\vbox{\vskip\ht#1\vskip\dp#1%
\pst@Verb{180 \tx@Rot}%
\vbox to \z@{\hbox to \z@{\box#1\hss}\vss}%
\pst@Verb{-180 \tx@Rot}}}}
\pslongbox{Rotateleft}{\psrotateleft}
\pslongbox{Rotateright}{\psrotateright}
\pslongbox{Rotatedown}{\psrotatedown}

\def\pst@starbox{%
\setbox\pst@hbox\hbox{\psframebox*[boxsep=false]{\unhbox\pst@hbox}}}
\def\pst@@makesmall#1{%
\setbox#1=\hbox to\z@{\hss\vbox to \z@{\vss\box#1\vss}\hss}}
\def\pst@@@makesmall#1{%
\pst@dimh=\psk@xref\wd#1%
\ifx\psk@yref\relax
\pst@dimg=\dp#1%
\else
\pst@dimg=\psk@yref\ht#1%
\advance\pst@dimg\psk@yref\dp#1%
\fi
\setbox#1=\hbox to\z@{%
\kern-\pst@dimh\vbox to\z@{\vss\box#1\kern-\pst@dimg}\hss}}
\def\psset@ref#1{\pst@expandafter\psset@@ref{#1}\@empty,,\@nil}
\def\psset@@ref#1#2,#3,#4\@nil{%
  \def\psk@xref{.5}%
  \def\psk@yref{.5}%
  \let\pst@makesmall\pst@@@makesmall
  \ifx\@empty#3\@empty
    \@nameuse{getref@#1}%
    \@nameuse{getref@#2}%
  \else
    \pst@checknum{#1#2}\psk@xref
    \pst@checknum{#3}\psk@yref
  \fi}
\def\getref@c{\let\pst@makesmall\pst@@makesmall}
\def\getref@t{\def\psk@yref{1}}
\def\getref@b{\def\psk@yref{0}}
\def\getref@B{\let\psk@yref\relax}
\def\getref@l{\def\psk@xref{0}}
\def\getref@r{\def\psk@xref{1}}
\psset@ref{c}
\def\psset@rot#1{%
\pst@expandafter{\@ifnextchar*{\psset@@@rot}{\psset@@rot}}{#1}\@nil}
\def\psset@@rot#1\@nil{%
\def\next##1@#1=##2@##3\@nil{%
\ifx\relax##2%
\pst@getangle{#1}\psk@rot
\else
\def\psk@rot{##2}%
\fi}%
\expandafter\next\pst@rottable @#1=\relax @\@nil}
\def\psset@@@rot#1#2\@nil{%
\psset@@rot#2\@nil
\edef\psk@rot{\pst@rotlist \ifx\psk@rot\@empty\else\psk@rot add \fi}}
\def\pst@rotlist{mark RAngle /a ED cleartomark a neg }
\def\pst@rottable{%
@0=%
@U=%
@L=90 %
@D=180 %
@R=-90 %
@N=\pst@rotlist
@W=\pst@rotlist 90 add %
@S=\pst@rotlist 180 add %
@E=\pst@rotlist 90 sub }
\psset@rot{0}
\def\tx@RotBegin{RotBegin }
\def\tx@RotEnd{RotEnd }
\def\pst@rotate#1#2{%
  \ifx#1\@empty\else
  \setbox#2=\hbox{\pst@Verb{#1 \tx@RotBegin}\box#2\pst@Verb{\tx@RotEnd}}%
  \fi%
}
\def\psput@cartesian#1{%
  \hbox to \z@{\kern\pst@dimg{\vbox to \z@{\vss\box#1\vskip\pst@dimh}\hss}}%
}
\def\psput@special#1{%
  \hbox{%
    \pst@Verb{{ \pst@coor } \tx@PutCoor \tx@PutBegin}%
    \box#1%
    \pst@Verb{\tx@PutEnd}%
  }%
}
\def\tx@PutCoor{PutCoor }
\def\tx@PutBegin{PutBegin }
\def\tx@PutEnd{PutEnd }
\def\rput{\def\pst@par{}\pst@ifstar{\@ifnextchar[{\rput@i}{\rput@ii}}}
\def\rput@i[#1]{\addto@par{ref={#1}}\rput@ii}
\def\rput@ii{\@ifnextchar({\rput@iv}{\rput@iii}}
\def\rput@iii#1{\addto@par{rot={#1}}\@ifnextchar({\rput@iv}{\rput@iv(\z@,\z@)}}
\def\rput@iv(#1){\pst@killglue\pst@makebox{\rput@v{#1}}}
\def\rput@v#1{%
  \begingroup
    \use@par
    \if@star\pst@starbox\fi
    \pst@makesmall\pst@hbox
    \pst@rotate\psk@rot\pst@hbox
    \psput@{#1}\pst@hbox
  \endgroup
  \ignorespaces%
}
\def\multirput{%
  \def\pst@par{}%
  \pst@ifstar{\@ifnextchar[{\multirput@i}{\multirput@ii}}%
}
\def\multirput@i[#1]{\addto@par{ref={#1}}\multirput@ii}
\def\multirput@ii{\@ifnextchar({\multirput@iv}{\multirput@iii}}
\def\multirput@iii#1{\addto@par{rot={#1}}\multirput@iv}
\def\multirput@iv(#1){%
  \@ifnextchar({\multirput@v(#1)}{\multirput@v(\z@,\z@)(#1)}%
}
\def\multirput@v(#1,#2)(#3,#4)#5{%
  \pst@makebox{\multirput@vi(#1,#2)(#3,#4){#5}}%
}
\def\multirput@vi(#1,#2)(#3,#4)#5{%
  \pst@killglue
  \begingroup
    \use@par
    \if@star\pst@starbox\fi
    \pst@makesmall\pst@hbox
    \pst@rotate\psk@rot\pst@hbox
    \pssetxlength\pst@dima{#1}%
    \pssetylength\pst@dimb{#2}%
    \pssetxlength\pst@dimc{#3}%
    \pssetylength\pst@dimd{#4}%
    \pst@cntg=#5\relax
    \pst@cnth=\@ne
    \leavevmode
    \loop
      \vbox to \z@{%
        \vss
        \hbox to \z@{\kern\pst@dima\copy\pst@hbox\hss}%
        \vskip\pst@dimb%
      }%
      \ifnum\pst@cntg>\pst@cnth
        \advance\pst@dima\pst@dimc
        \advance\pst@dimb\pst@dimd
        \advance\pst@cnth\@ne
    \repeat 
    \endgroup
  \ignorespaces%
}
\newif\if@fixedradius
\def\cput{\pst@object{cput}}
\def\cput@i{\@fixedradiusfalse\cput@ii}
\def\cput@ii{\pst@killglue\@ifnextchar({\cput@iv}{\cput@iii}}
\def\cput@iii#1{%
  \addto@par{rot={#1}}%
  \@ifnextchar({\cput@iv}{\cput@iv(\z@,\z@)}%
}
\def\cput@iv(#1){\pst@makebox{\cput@v{#1}}}
  \def\cput@v#1{%
  \begingroup
    \use@par
    \setbox\pst@hbox=\hbox{%
      \psboxsepfalse
      \if@fixedradius\psCirclebox@ii\else\pscirclebox@ii\fi%
    }%
    \pst@@makesmall\pst@hbox
    \pst@rotate\psk@rot\pst@hbox
    \psput@{#1}\pst@hbox
  \endgroup
  \ignorespaces%
}
\def\Cput{\pst@object{Cput}}
\def\Cput@i{\@fixedradiustrue\cput@ii}
\newdimen\pslabelsep
\def\psset@labelsep#1{\pssetlength\pslabelsep{#1}}
\psset@labelsep{5pt}
\def\psset@refangle#1{\pst@expandafter\psset@@refangle{#1}\@nil}
\def\psset@@refangle#1\@nil{%
\def\next##1@#1=##2"##3@##4\@nil{%
\ifx\relax##2%
\pst@getangle{#1}\psk@refangle
\def\psk@uputref{}%
\else
\def\psk@refangle{##2 }%
\def\psk@uputref{##3}%
\fi}%
\expandafter\next\pst@refangletable @#1=\relax"@\@nil}
\def\pst@refangletable{%
@r=0"20%
@u=90"02%
@l=180"10%
@d=-90"01%
@ur=45"22%
@ul=135"12%
@dr=-135"21%
@dl=-45"11}
\psset@refangle{0}
\def\uput{\def\pst@par{}\pst@ifstar{\@ifnextchar[{\uput@ii}{\uput@i}}} 
\def\uput@i#1{\addto@par{labelsep=#1}\uput@ii}
\def\uput@ii[#1]{%
\addto@par{refangle={#1}}%
\@ifnextchar({\uput@iv}{\uput@iii}}
\def\uput@iii#1{%
\addto@par{rot={#1}}%
\@ifnextchar({\uput@iv}{\uput@iv(\z@,\z@)}}
\def\uput@iv(#1){\pst@killglue\pst@makebox{\uput@v{#1}}}
\def\uput@v#1{%
\begingroup
\use@par
\if@star\pst@starbox\fi
\uput@vi
\psput@{#1}\pst@hbox
\endgroup
\ignorespaces}
\def\uput@vi{%
\ifx\psk@uputref\@empty
\uput@vii\tx@UUput{}%
\else
\ifx\psk@rot\@empty
\expandafter\uput@viii\psk@uputref
\else
\uput@vii\tx@UUput{}%
\fi
\fi}
\def\uput@vii#1#2{%
  \edef\pst@coor{%
    \pst@number\pslabelsep
    #2%
    \pst@number{\wd\pst@hbox}%
    \pst@number{\ht\pst@hbox}%
    \pst@number{\dp\pst@hbox}%
    \psk@refangle\space \ifx\psk@rot\@empty\else\psk@rot\space sub \fi
    \tx@Uput #1}%
  \setbox\pst@hbox=\hbox to\z@{\hss\vbox to\z@{\vss\box\pst@hbox\vss}\hss}%
  \setbox\pst@hbox=\psput@special\pst@hbox
  \ifx\psk@rot\@empty\else\pst@rotate\psk@rot\pst@hbox\fi}
\def\uput@viii#1#2{%
  \ifnum#1>\z@\ifnum#2>\z@\pslabelsep=.707\pslabelsep\fi\fi
  \setbox\pst@hbox=\vbox to\z@{%
    \ifnum#2=1 \vskip\pslabelsep\else\vss\fi
    \hbox to\z@{%
      \ifnum#1=2 \hskip\pslabelsep\else\hss\fi
      \box\pst@hbox
      \ifnum#1=1 \hskip\pslabelsep\else\hss\fi}%
    \ifnum#2=2 \vskip\pslabelsep\else\vss\fi}}
\def\tx@Uput{Uput }
\def\tx@UUput{UUput }
\def\Rput{\def\pst@par{}\pst@ifstar{\@ifnextchar[{\Rput@ii}{\Rput@i}}}
\def\Rput@i#1{\addto@par{labelsep=#1}\Rput@ii}
\def\Rput@ii[#1]{\addto@par{ref={#1}}\@ifnextchar({\Rput@iv}{\Rput@iii}}
\def\Rput@iii#1{\addto@par{rot={#1}}\@ifnextchar({\Rput@iv}{\Rput@iv(\z@,\z@)}}
\def\Rput@iv(#1){\pst@killglue\pst@makebox{\Rput@v{#1}}}
\def\Rput@v#1{%
\begingroup
\use@par
\if@star\pst@starbox\fi
\Rput@vi
\pst@makesmall\pst@hbox
\pst@rotate\psk@rot\pst@hbox
\psput@{#1}\pst@hbox
\endgroup
\ignorespaces}
\def\Rput@vi{%
\pst@dimg=\dp\pst@hbox
\advance\pst@dimg\pslabelsep
\dp\pst@hbox=\pst@dimg
\pst@dimg=\ht\pst@hbox
\advance\pst@dimg\pslabelsep
\ht\pst@hbox=\pst@dimg
\setbox\pst@hbox\hbox{\kern\pslabelsep\box\pst@hbox\kern\pslabelsep}}%
\def\oldpsput{%
\def\pst@par{}\pst@ifstar{\@ifnextchar[{\oldpsput@i}{\oldpsput@ii}}}
\def\oldpsput@i[#1]{\addto@par{ref={#1}}\oldpsput@ii}
\def\oldpsput@ii{\@ifnextchar<{\oldpsput@iii}{\oldpsput@iv}}
\def\oldpsput@iii<#1>{\rput@iii{#1}}
\def\OldPsput{\let\psput\oldpsput}
\def\NewPsput{\let\psput\rput}
%
\newpsstyle{gridstyle}{subgriddiv=0,gridcolor=lightgray,griddots=10,gridlabels=8pt}
\newif\ifshowgrid
\def\psset@showgrid#1{\@nameuse{showgrid#1}}
\psset@showgrid{false}
\newdimen\psk@shift
\def\psset@shift#1{\pssetlength\pst@dimg{#1}%
  \psk@shift\pst@dimg}
\psset@shift{0}
%
\def\pspicture{\begingroup\pst@ifstar\pst@picture}
\def\pst@picture{%
\@ifnextchar[{\pst@@picture}{\pst@@picture[]}}
\def\pst@@picture[#1]#2(#3,#4){%
\@ifnextchar({\pst@@@picture[#1](#3,#4)}%
{\pst@@@picture[#1](0,0)(#3,#4)}}
\def\pst@@@picture[#1](#2,#3)(#4,#5){%
  \pssetxlength\pst@dima{#2}%
  \pssetylength\pst@dimb{#3}%
  \pssetxlength\pst@dimc{#4}%
  \pssetylength\pst@dimd{#5}%
  \ifdim\pst@dima>\pst@dimc%
    \pst@dimg=\pst@dima%
    \pst@dima=\pst@dimc%
    \pst@dimc=\pst@dimg%
  \fi%
  \ifdim\pst@dimb>\pst@dimd%
    \pst@dimg=\pst@dimb%
    \pst@dimb=\pst@dimd%
    \pst@dimd=\pst@dimg%
  \fi%
  \setbox\pst@hbox=\hbox\bgroup%
  \begingroup\KillGlue%
  \@ifundefined{@latexerr}{}{\let\unitlength\psunit}%
  \edef\pic@coor{(#2,#3)(#2,#3)(#4,#5)}%
  \psset{showgrid=false}
  \def\pst@tempA{#1}%
  \ifx\pst@tempA\@empty\else\psset{#1}\fi
  \ifshowgrid\psgrid[style=gridstyle]\fi%
}
\def\pic@coor{(0,0)(0,0)(10,10)}
\newdimen\pst@shift
\def\endpspicture{%
  \pst@killglue
  \global\pst@shift=\psk@shift
  \endgroup
  \egroup
  \ifdim\wd\pst@hbox=\z@\else
  \fi
  \ht\pst@hbox=\pst@dimd
  \dp\pst@hbox=-\pst@dimb
  \setbox\pst@hbox=\hbox{%
    \kern-\pst@dima
    \pst@dimd-\pst@shift
    \advance\pst@dimd\pst@dimb
    \lower\pst@dimd%
    \box\pst@hbox%
    \kern\pst@dimc}%
  \if@star\setbox\pst@hbox=\hbox{\clipbox@@\z@}\fi
  \leavevmode\box\pst@hbox
  \endgroup%
  \global\psk@shift\z@
}
\@namedef{pspicture*}{\pspicture*}
\@namedef{endpspicture*}{\endpspicture}
\def\tx@BeginOL{BeginOL }
\def\tx@InitOL{InitOL }
\def\pst@initoverlay#1{\pst@Verb{\tx@InitOL /TheOL (#1) def}}
\def\AltOverlayMode{%
  \def\pst@initoverlay##1{%
    \pst@Verb{%
      \tx@InitOL
      /Visible { initclip } def
      /Invisible {
        CP newpath OLUnit itransform moveto clip newpath moveto
      } def
      /TheOL (##1) def}}}
\def\pst@overlay#1{%
  \edef\curr@overlay{#1}%
  \pst@Verb{(#1) BOL}%
  \aftergroup\pst@endoverlay}
\def\pst@endoverlay{%
  \pst@Verb{(\curr@overlay) BOL}}
\def\curr@overlay{all}
\newbox\theoverlaybox
\def\overlaybox{%
  \global\setbox\theoverlaybox=\hbox\bgroup
  \begingroup
  \let\psoverlay\pst@overlay
  \def\overlaybox{%
    \@pstrickserr{Overlays cannot be nested}\@eha}%
  \def\putoverlaybox{%
    \@pstrickserr{You must end the overlay box
         before using \string\putoverlaybox}}%
  \psoverlay{main}%
  \ignorespaces}
\def\endoverlaybox{\endgroup\egroup}
\def\putoverlaybox#1{%
\hbox{\pst@initoverlay{#1}\copy\theoverlaybox}}
\def\psoverlay{\@pstrickserr{\string\psoverlay\space
can only be used after \string\overlaybox}}
\ifx\pstcustomize\relax \input pstricks.con \fi
\catcode`\@=\PstAtCode\relax
%
 
\input xy \xyoption{all}
\def\Pic{{\rm Pic}}
\newwrite\num
%
\output={\if N\header\headline={\hfill}\fi
\plainoutput\global\let\header=Y}
\magnification\magstep1
\tolerance = 500
\hsize=14.4true cm
\vsize=22.5true cm
\parindent=6true mm\overfullrule=2pt
\newcount\kapnum \kapnum=0
\newcount\parnum \parnum=0
\newcount\procnum \procnum=0
\newcount\nicknum \nicknum=1
\font\ninett=cmtt9

\font\ninebf=cmbx9

\font\sixbf=cmbx6
\font\ninesl=cmsl9

\font\nineit=cmti9

\font\ninerm=cmr9

\font\sixrm=cmr6
\font\ninei=cmmi9
\font\eighti=cmmi8
\font\sixi=cmmi6
\skewchar\ninei='177 \skewchar\eighti='177 \skewchar\sixi='177
\font\ninesy=cmsy9
\font\eightsy=cmsy8
\font\sixsy=cmsy6
\skewchar\ninesy='60 \skewchar\eightsy='60 \skewchar\sixsy='60
\font\titelfont=cmr10 scaled 1440
\font\paragratit=cmbx10 scaled 1200

\font\name=cmcsc10
\font\emph=cmbxti10

\font\tenmsbm=msbm10
\font\sevenmsbm=msbm7
%

%

%
\font\teneufm=eufm10
\font\seveneufm=eufm7
\font\fiveeufm=eufm5
\newfam\eufmfam
\textfont\eufmfam=\teneufm
\scriptfont\eufmfam=\seveneufm
\scriptscriptfont\eufmfam=\fiveeufm

\font\tenmsam=msam10
\font\sevenmsam=msam7
\font\fivemsam=msam5
\newfam\msamfam
\textfont\msamfam=\tenmsam
\scriptfont\msamfam=\sevenmsam
\scriptscriptfont\msamfam=\fivemsam
\font\tenmsbm=msbm10
\font\sevenmsbm=msbm7
\font\fivemsbm=msbm5
\newfam\msbmfam
\textfont\msbmfam=\tenmsbm
\scriptfont\msbmfam=\sevenmsbm
\scriptscriptfont\msbmfam=\fivemsbm
\def\Bbb#1{{\fam\msbmfam\relax#1}}
\def\cz{{\kern0.4pt\Bbb C\kern0.7pt}
}
\def\ez{{\kern0.4pt\Bbb E\kern0.7pt}
}
\def\fz{{\kern0.4pt\Bbb F\kern0.3pt}}
\def\gz{{\kern0.4pt\Bbb Z\kern0.7pt}}
\def\hz{{\kern0.4pt\Bbb H\kern0.7pt}
}
\def\kz{{\kern0.4pt\Bbb K\kern0.7pt}
}
\def\nz{{\kern0.4pt\Bbb N\kern0.7pt}
}
\def\oz{{\kern0.4pt\Bbb O\kern0.7pt}
}
\def\rz{{\kern0.4pt\Bbb R\kern0.7pt}
}
\def\sz{{\kern0.4pt\Bbb S\kern0.7pt}
}
\def\pz{{\kern0.4pt\Bbb P\kern0.7pt}
}
\def\qz{{\kern0.4pt\Bbb Q\kern0.7pt}
}
\newskip\ttglue
\def\ninepoint{\def\rm{\fam0\ninerm}%
  \textfont0=\ninerm \scriptfont0=\sixrm \scriptscriptfont0=\fiverm
  \textfont1=\ninei \scriptfont1=\sixi \scriptscriptfont1=\fivei
  \textfont2=\ninesy \scriptfont2=\sixsy \scriptscriptfont2=\fivesy
  \textfont3=\tenex \scriptfont3=\tenex \scriptscriptfont3=\tenex
  \def\it{\fam\itfam\nineit}%
  \textfont\itfam=\nineit
  \def\sl{\fam\slfam\ninesl}%
  \textfont\slfam=\ninesl
  \def\bf{\fam\bffam\ninebf}%
  \textfont\bffam=\ninebf \scriptfont\bffam=\sixbf
   \scriptscriptfont\bffam=\fivebf
  \def\tt{\fam\ttfam\ninett}%
  \textfont\ttfam=\ninett
  \tt \ttglue=.5em plus.25em minus.15em
  \normalbaselineskip=11pt
  \font\name=cmcsc9
  \let\sc=\sevenrm
  \let\big=\ninebig
  \setbox\strutbox=\hbox{\vrule height8pt depth3pt width0pt}%
  \normalbaselines\rm
  \def\sl{\it}}

\headline={\ifodd\pageno\rightheadline\else\leftheadline\fi}
\def\rightheadline{\ninepoint Paragraphen"uberschrift\hfill\folio}
\def\leftheadline{\ninepoint\folio\hfill Chapter"uberschrift}
\let\header=Y
\def\titel#1{\need 9cm \vskip 2truecm
\parnum=0\global\advance \kapnum by 1
{\baselineskip=16pt\lineskip=16pt\rightskip0pt
plus4em\spaceskip.3333em\xspaceskip.5em\pretolerance=10000\noindent
\titelfont Chapter \uppercase\expandafter{\romannumeral\kapnum}.
#1\vskip2true cm}\def\leftheadline{\ninepoint
\folio\hfill Chapter \uppercase\expandafter{\romannumeral\kapnum}.
#1}\let\header=N
}
\def\Titel#1{\need 9cm \vskip 2truecm
\global\advance \kapnum by 1
{\baselineskip=16pt\lineskip=16pt\rightskip0pt
plus4em\spaceskip.3333em\xspaceskip.5em\pretolerance=10000\noindent
\titelfont\uppercase\expandafter{\romannumeral\kapnum}.
#1\vskip2true cm}\def\leftheadline{\ninepoint
\folio\hfill\uppercase\expandafter{\romannumeral\kapnum}.
#1}\let\header=N
}
\def\need#1cm {\par\dimen0=\pagetotal\ifdim\dimen0<\vsize
\global\advance\dimen0by#1 true cm
\ifdim\dimen0>\vsize\vfil\eject\noindent\fi\fi}
\def\neupara#1{\par\penalty-2000
\procnum=0\global\advance\parnum by 1
\vskip1cm\noindent{\paragratit \the\parnum. #1}%
\def\rightheadline{\ninepoint\S\the\parnum.\ #1\hfill \folio}%
\vskip 8mm\noindent}
\def\Proclaim #1 #2\finishproclaim {\bigbreak\noindent
{\bf#1\unskip{}. }{\it#2}\medbreak\noindent}
%
\gdef\proclaim #1 #2 #3\finishproclaim {\bigbreak\noindent%
\global\advance\procnum by 1
{%
{\relax\ifodd \nicknum
\hbox to 0pt{\vrule depth 0pt height0pt width\hsize
   \quad \ninett#3\hss}\else {}\fi}%
\bf\the\parnum.\the\procnum\ #1\unskip{}. }
{\it#2}
\immediate\write\num{\string\def
 \expandafter\string\csname#3\endcsname
 {\the\parnum.\the\procnum}}
\medbreak\noindent}
\newcount\stunde \newcount\minute \newcount\hilfsvar
\def\uhrzeit{
    \stunde=\the\time \divide \stunde by 60
    \minute=\the\time
    \hilfsvar=\stunde \multiply \hilfsvar by 60
    \advance \minute by -\hilfsvar
    \ifnum\the\stunde<10
    \ifnum\the\minute<10
    0\the\stunde:0\the\minute~Uhr
    \else
    0\the\stunde:\the\minute~Uhr
    \fi
    \else
    \ifnum\the\minute<10
    \the\stunde:0\the\minute~Uhr
    \else
    \the\stunde:\the\minute~Uhr
    \fi
    \fi
    }
\def\calA{{\cal A}} \def\calB{{\cal B}}

\def\calC{{\cal C}} \def\calD{{\cal D}}
 
\def\calG{{\cal G}}

 \def\calX{{\cal X}}
\def\calY{{\cal Y}} \def\calZ{{\cal Z}}

\def\dim{\mathop{\rm dim}\nolimits}

\def\Pic{\mathop{\rm Pic}\nolimits}

\def\boxit#1{
  \vbox{\hrule\hbox{\vrule\kern6pt
  \vbox{\kern8pt#1\kern8pt}\kern6pt\vrule}\hrule}}
\def\Boxit#1{
  \vbox{\hrule\hbox{\vrule\kern2pt
  \vbox{\kern2pt#1\kern2pt}\kern2pt\vrule}\hrule}}
\def\rahmen#1{
  $$\boxit{\vbox{\hbox{$ \displaystyle #1 $}}}$$}

\def\smallni{\smallskip\noindent }
\def\medni{\medskip\noindent }

\def\loma{\longmapsto}
\def\imag{{\rm i}}

\def\square{\hbox{\hbox to 0pt{$\sqcup$\hss}\hbox{$\sqcap$}}}
\def\qed{\ifmmode\square\else{\unskip\nobreak\hfil
\penalty50\hskip3em\null\nobreak\hfil\square
\parfillskip=0pt\finalhyphendemerits=0\endgraf}\fi}
\def\pn{\the\parnum.\the\procnum}
\def\downmapsto{{\buildrel
        {\vbox{\hbox{\hskip.2pt$\scriptstyle-$}}}
        \over{\raise7pt\vbox{\vskip-4pt\hbox{$\textstyle\downarrow$}}}}}

\input calabi-phys.num
\nopagenumbers
\immediate\newwrite\num
\nicknum=0  
\let\header=N

\immediate\openout\num=calabi-phys.num
\immediate\newwrite\num\immediate\openout\num=calabi-phys.num
\def\RAND#1{\vskip0pt\hbox to 0mm{\hss\vtop to 0pt{%
  \raggedright\ninepoint\parindent=0pt%
  \baselineskip=1pt\hsize=2cm #1\vss}}\noindent}
\noindent
\centerline{\titelfont A rigid Calabi-Yau manifold with Picard number two}%
\def\leftheadline{\ninepoint\folio\hfill
A rigid Calabi-Yau manifold with Picard number two}%
\def\rightheadline{\ninepoint Introduction\hfill \folio}%
\headline={\ifodd\pageno\rightheadline\else\leftheadline\fi}
\vskip 5mm
\centerline{Eberhard Freitag,   2015}%
\vskip5mm\noindent%
\let\header=N%
\def\imag{{\rm i}}%
\def\transpose#1{\kern1pt{^t\kern-3pt#1}}%
\def\t{\vartheta}%
\centerline{\vbox{\noindent\hsize=10cm
{\ninepoint{\bf Abstract}
\smallni
We study a projective Calabi--Yau threefold $\tilde\calY$
which has been constructed in [FS]. It is rigid ($h^{12}=0$) and has Picard
number $h^{11}=2$. We construct a pair of divisors $\calD^\pm$ which give a basis of
$\Pic(\tilde\calY)\otimes_\gz\qz$
and determine all intersection numbers $\calD^\pm\cdot\calD^\pm\cdot\calD^\pm$.}}}
\vskip5mm\noindent
We want to thank Sergey V.~Ketov, who motivated us to determine the intersection numbers
which play a role for physical applications in string theory. We are  grateful to
S.~Cynk who explained us several details about the resolution of nodes and we thank also
R.~Salvati Manni for widespread discussions and suggestions.
\vskip5mm\noindent
{\paragratit Introduction}%
\medni
Basic for our example is a certain complete intersection $\calX$
of four quadrics introduced in
the paper [GN] of van Geemen and Nygaard:
$$\eqalign{Y_0^2&=X_0^2+X_1^2+X_2^2+X_3^2,\cr
Y_1^2&=X_0^2-X_1^2+X_2^2-X_3^2,\cr
Y_2^2&=X_0^2+X_1^2-X_2^2-X_3^2,\cr
Y_3^2&=X_0^2-X_1^2-X_2^2+X_3^2.\cr}\leqno\qquad\calX:$$
The variety $\calX$ has 96 isolated singularities which are ordinary
double points (nodes). 
\smallskip
In the paper [CM] it has
been pointed out that the results of [GN] imply that
$\calX$ admits a resolution that is a (projective)
Calabi--Yau threefold. 
The basic result -- essentially due to van Geemen and
Nygaard [GN] -- is the following theorem.
\Proclaim
{Theorem}
{The Hodge numbers of a Calabi--Yau desingularization of $\calX$ are
$$h_{11}=32,\qquad h_{12}=0.$$
Hence this Calabi--Yau manifold is rigid.}
\finishproclaim
In the paper [FS], we constructed a certain group $G$ of order 16 of biholomorphic
automorphisms of $\calX$ that acts freely on $\calX$.
The basic thing is that there exists a {\it projective\/} 
small resolution $\tilde\calX$ such that $G$ extends
as group of biholomorphic mappings 
to $\tilde\calX$.
The quotient
$\tilde\calY=\tilde\calX/G$ then is a projective resolution  
of $\calY:=\calX/G$ in form of a rigid Calabi-Yau manifold
whose Picard number is two. The varieties $\calX$ and $\calY$ are Siegel modular
three folds. But for this paper the modular background is not necessary. Everything
can be done by using the explicit equations.
\smallskip
We recall the definition of the group $G$. Let $P_1,\dots,P_8$ be homogenous polynomials
in $\cz[Y_0,\dots,Y_3,X_0,\dots,X_3]$ of the same degree such that not all of them vanish.
Then we can consider
the rational map from $P^7$ into itself,
$$(y_0,\dots,y_3,x_0,\dots,x_3)\loma (P_1(y_0,\dots,x_3),\dots,P_8(y_0,\dots,x_3)).$$
We denote this map symbolically by
$$(P_1,\dots,P_8).$$
\Proclaim
{Theorem}
{The two transformations
$$\eqalign{&(-Y_2,\imag Y_3,\imag Y_0,Y_1,-\imag X_3,-\imag X_2,X_1,-X_0),\cr&
(\imag Y_1,Y_0,-\imag Y_3,Y_2,X_1,\imag X_0,X_3,-\imag X_2)}$$
define biholomorphic transformations of $\calX$ onto itself.
The group $G$ which is generated by the two biholomorphic transformations
is isomorphic to $\gz/4\times\gz/4$ and  acts freely on $\calX$. The quotient $\calY=\calX/G$
has a resolution $\tilde\calY$ in the form of
a (projective) rigid Calabi--Yau manifold  ($h^{12}=0$) with Euler number $e=4$
and Picard number $h^{11}=2$.}
\finishproclaim
{\it Proof.} The proof is indicated in [FS]. 
We give some more details.
The  matrix group $\calG$, generated by the two matrices
$$\pmatrix{
0&0&-1&0& 0&0&0&0\cr
0&0&0&\imag& 0&0&0&0\cr
\imag&0&0&0& 0&0&0&0\cr
0&1&0&0& 0&0&0&0\cr
0&0&0&0& 0&0&0&-\imag\cr
0&0&0&0& 0&0&-\imag&0\cr
0&0&0&0& 0&1&0&0\cr
0&0&0&0& -1&0&0&0\cr
}
\qquad
\pmatrix{
0&\imag&0&0& 0&0&0&0\cr
1&0&0&0& 0&0&0&0\cr
0&0&0&-\imag& 0&0&0&0\cr
0&0&1&0& 0&0&0&0\cr
0&0&0&0& 0&1&0&0\cr
0&0&0&0& \imag&0&0&0\cr
0&0&0&0& 0&0&0&1\cr
0&0&0&0& 0&0&-\imag&0\cr
},$$
has order 64. Its center $\calC$ is the group of order 4,
generated by the matrix $\imag E$ (where $E$ denotes the unit matrix).
There is a natural homomorphism $\calG\to G$ which gives the identification
$\calG/\calZ\cong G$.
(As usual, $n\times n$-matrices act from the left on $\cz^n$: the elements of
$\cz^n$ are written as columns and then matrices act by matrix multiplication from
the left).
\smallskip
Next we prove that $G$ acts fixed point free on $\calX$. It is enough to do
this for elements of order two in $G$. There are three of them, namely the two squares
of the two generators and the product of the two squares. All three are
diagonal matrices with two different entries. It is very easy to check that they
have no fixed point in $\calX$.
\qed
\neupara{Intersection numbers between basic divisors}%
In the paper [FS], 188 two-dimensional  subvarieties of $\calX$
have been defined. We use here two of them.
They are components of the hyperplane section $X_0+\cdots+X_3=0$.
One of them is cut out by the equations
$$D^+:\quad X_0+X_1+X_2+X_3=Y_1(X_1+X_3)+(\sqrt2/2) Y_2Y_3=0.$$
The vanishing ideal of this variety is the radical of the ideal generated by
the defining ideals of $\calX$ and $D^+$. It can be computed as
$$\eqalign{
    &Y_0^2 - 2X_1^2 - 2X_1X_2 - 2X_1X_3 - 2X_2^2 - 2X_2X_3 - 2X_3^2,\cr
    &Y_1^2 - 2X_1X_2 - 2X_1X_3 - 2X_2^2 - 2X_2X_3,\cr
    &Y_1Y_2 + \sqrt2 Y_3X_1 + \sqrt2 Y_3X_2,\cr
    &Y_1Y_3 + \sqrt2 Y_2X_2 + \sqrt2 Y_2X_3,\cr
    &Y_1X_1 + Y_1X_3 + 1/2\sqrt2 Y_2Y_3,\cr
    &Y_2^2 - 2X_1^2 - 2X_1X_2 - 2X_1X_3 - 2X_2X_3,\cr
    &Y_3^2 - 2X_1X_2 - 2X_1X_3 - 2X_2X_3 - 2X_3^2,\cr
    &X_0 + X_1 + X_2 + X_3.\cr}$$
Replacing $\sqrt2$ by $-\sqrt2$ we get a complementary ideal
$$D^-:\quad X_0+X_1+X_2+X_3=Y_1(X_1+X_3)-(\sqrt2/2) Y_2Y_3=0.$$
In  the following we use the notation $D^\pm$ for the corresponding irreducible
divisors.
We recall some basics about of divisors and their intersection numbers.
\smallskip
Let $X$ be a normal  irreducible algebraic variety of dimension $n$ over $\cz$.
A divisor on $X$ is a formal sum of irreducible closed subvarieties of codimension one.
Let $Y\subset X$ be a subvariety of everywhere
codimension $1$. The associated divisor is the sum of all irreducible components
(with multiplicities $1$). 
To every rational function its principal divisor can be associated.
Two divisors are called equivalent if their difference is principal.
Since the singular locus of $X$ has codimension $\ge2$, the divisors 
(divisor classes) on $X$ 
are in one-to-one
correspondence with the divisors (divisor classes) on the regular locus. 
Let $\pi:\tilde X\to X$ be a 
small resolution. Small means that there exists a closed subvariety $T\subset\tilde X$
of codimension $\ge 2$ such that $\pi$ defines a biholomorphic map
of $\tilde X-T$ onto the regular locus of $X$. The divisors (divisor classes) of $X$ are in
one-to-one correspondence with those of $\tilde X$.
\smallskip
Now we assume that $X$ is projective and non-singular.
Then the intersection number $D_1\cdots D_n$ 
of $n$ divisors $D_1,\dots D_n$ 
can be defined. It is invariant under 
equivalence and it is $\gz$-multilinear. In the case that the divisors are effective
and intersect only in finitely many points $P_1,\dots, P_m$,
the intersection number is the sum of the multiplicities of $P_i$ in the
scheme theoretic intersection $D_1\cap\dots\cap D_n$. For details
we refer to [Sh]. 
\smallskip
We consider a projective small resolution $\tilde\calX$. (This is a Calabi--Yau manifold).
We can consider divisors $D_1,D_2,D_3$ on $\calX$ as divisors on $\tilde\calX$ 
and then study intersection
numbers of three divisors. 
To be precise, $D_1\cdot D_2\cdot D_3$ means the intersection number on $\tilde\calX$.
These numbers may depend on the choice of the resolution.
But, if one of the three divisors is equivalent to a divisor which 
avoids all nodes, then the intersection
is independent of the choice. For example, this is the case for hyper plane sections. 
They are all equivalent
and  there is one which avoids a finite number of given points.
\smallskip
We consider the divisors $D^+$ and $D^-$ on $\calX$. Their sum $H$ is the divisor
of a hyperplane section. We can consider the three divisors also as divisors
on $\tilde\calX$.
In this section we want to compute the intersection numbers of three of the divisors
$H,D^+,D^-$ (on $\tilde\calX$). As we explained, they to not depend on the choice
of the small resolution.
\proclaim
{Proposition}
{Let $H$ be a divisor on $\calX$ that represents a hyperplane section, i.e.\ $H\sim (X_0)$.
Then
$$H\cdot H\cdot H=16.$$
}
HHH%
\finishproclaim
{\it Proof.\/} The hyperplane sections $(X_0)$, $(X_1)$, $(X_2)$ have (on $\calX$) 16
intersection points,
$$[\pm1,\pm\imag,\pm\imag,\pm1,0,0,0,1].$$
None of them is a node. Their multiplicity is one.\qed
\proclaim
{Proposition}
{Let $H$ be a divisor on $\calX$ that represents a hyperplane section, i.e.\ $H\sim (X_0)$.
Then
$$H\cdot H\cdot D^+=H\cdot H\cdot D^-=8.$$
}
HHDplus%
\finishproclaim
{\it Proof.\/}
There  is a biholomorphic transformation of $\calX$ that sends $Y_2$ to $-Y_2$ and preserves 
the remaining
variables. The divisor class of $H$ remains fixed and $D^+$ and $D^-$ are interchanged. Hence
$H\cdot H\cdot D^+=H\cdot H\cdot D^-$. Their sum is $16$ by Proposition \HHH.
(We should mention that a biholomorphic transformation usually does not extend to $\tilde\calX$.
But there is a second resolution $\bar\calX$ such that it extends to a biholomorphic 
map $\tilde X\to\bar\calX$. This shows that $H\cdot H\cdot D^+$, computed on $\tilde\calX$,
equals $H\cdot H\cdot D^-$, computed on $\bar\calX$. But as we noted already these intersection
numbers are independent of the choice of a small resolution.)
\qed
\proclaim
{Proposition}
{Let $H$ be a hyperplane section, i.e.\ $H\sim (X_0)$.
Then
$$H\cdot D^+\cdot D^-=20,\quad H\cdot D^+\cdot D^+=H\cdot D^-\cdot D^-=-12.$$
}
HDD%
\finishproclaim
{\it Proof.\/}
The divisors $(Y_0)$, $D^+$, $D^-$ have 12
intersection points, None of them is a node. The 4 points
$$[0,\pm2\imag,0,0,-\imag,-1,\imag,1],\quad [0,\pm2\imag,0,0,\imag,-1,-\imag,1]$$
have multiplicity 3. The remaining 8
$$\eqalign{&
[0,0,\pm2\imag,0,-\imag,\imag,-1,1],\quad [0,0,\pm2\imag,0,\imag,-\imag,-1,1],\cr&
[0,0,0,\pm2,-1,\imag,-\imag,1],\quad[0,0,0,\pm2,-1,-\imag,\imag,1]\cr}
$$
have multiplicity $1$. This gives the intersection number $4\cdot 3+8=20$.
The rest  follows with the help of Proposition \HHDplus.
\qed
\proclaim
{Proposition}
{We have
$$\eqalign{D^+\cdot D^+\cdot D^+=&D^-\cdot D^-\cdot D^-=-22,\cr
\quad D^+\cdot D^+\cdot D^-=&D^+\cdot D^-\cdot D^-=10.\cr}$$
}
DDD%
\finishproclaim
{\it Proof.\/} We have
$$D^+\cdot D^+\cdot D^-=D^+\cdot D^-\cdot D^-.$$
The sum of both is 20 (Proposition \HDD). Hence we obtain the value 10.
We also have
$$D^+\cdot D^+\cdot D^+=D^+\cdot\ D^+\cdot H-D^+\cdot D^+\cdot\ D^-=-22.$$
This completes the proof of Proposition \DDD.\qed
\neupara{Intersection numbers between translates}%
\proclaim
{Proposition}
{We have
$$H\cdot H\cdot g(D^+)=H\cdot H\cdot g(D^-)=8.$$}
HHg%
\finishproclaim
This follows from Proposition \HHDplus.\qed
\proclaim
{Proposition}
{We have
$$H\cdot g(D^+)\cdot g(D^-)=20,\quad  H\cdot g(D^+)\cdot g(D^+)=
H\cdot g(D^-)\cdot g(D^-)=-12.$$}
Hgga%
\finishproclaim
This follows from Proposition \HDD.\qed
\proclaim
{Proposition}
{We have
$$H\cdot g_1(D^+)\cdot g_2(D^+) =H\cdot g_1(D^+)\cdot g_2(D^-) =H\cdot g_1(D^-)\cdot g_2(D^-) =4$$
for all different $g_1\ne g_2$ in $G$.}
Hgg%
\finishproclaim
{\it Proof.\/} It is sufficient to treat the case where $g_1$ is the identity. Then
$g=g_2$ is different from the identity. In this case the intersection
of $Y_0=0$, $D^+$ and $g(D^+)$ consists of 4 points.
They are members of the following list of 12 points.
\hbox{\vbox{\hsize 3.8cm
$$\eqalign{
&[0,0,0,-2,-1,-\imag,\imag],\cr
&[0,0,2\imag,0,\imag,-\imag,-1],\cr
&[0,0,2\imag,0,-\imag,\imag,-1],\cr
&[0,2\imag,0,0,-\imag,-1,\imag],\cr}$$}
\vbox{\hsize=3.8cm
$$\eqalign{
&[0,0,-2\imag,0,\imag,-\imag,-1],\cr
&[0,2\imag,0,0,\imag,-1,-\imag],\cr
&[0,-2\imag,0,0,-\imag,-1,\imag],\cr
&[0,0,0,-2,-1,\imag,-\imag],\cr}$$}
\vbox{\hsize=3.8cm
$$\eqalign{
&[0,0,-2\imag,0,-\imag,\imag,-1],\cr
&[0,0,0,2,-1,-\imag,\imag],\cr
&[0,0,0,2,-1,\imag,-\imag],\cr
&[0,-2\imag,0,0,\imag,-1,-\imag].\cr}$$}}
\smallni
None of them is a a node.
They have multiplicity one.
This shows $H\cdot g_1(D^+)\cdot g_2(D^+)=4$. The rest follows from Proposition \HHg.
\qed
\proclaim
{Proposition}
{For $g\in G$ we have
$$\eqalign{g(D^+)\cdot g(D^-)\cdot g(D^-)=&g(D^+)\cdot g(D^+)\cdot g(D^-)=10,\cr
g(D^+)\cdot g(D^+)\cdot g(D^+)=&g(D^-)\cdot g(D^-)\cdot g(D^-)=-22.\cr}$$
}
gDDD%
\finishproclaim
This follows from Proposition \DDD.\qed
\proclaim
{Proposition}
{Let $g_1,g_2,g_3$ be three  elements from $G$ which contain two different
ones, $\#\{g_1,g_2,g_3\}=2$. Then we have
$$g_1(D^\pm)\cdot g_2(D^\pm)\cdot g_3(D^\pm)=10$$
if the three factors $g_i(D^\pm)$ are pairwise different. For all other
sign combinations we have 
$$g_1(D^\pm)\cdot g_2(D^\pm)\cdot g_3(D^\pm)=-6.$$
}
gegz%
\finishproclaim
{\it Proof.\/} Applying the transformation $Y_2\mapsto -Y_2$, we see
$g_1(D^+)\cdot g_2(D^+)\cdot g_2(D^-)=g_1(D^-)\cdot g_2(D^+)\cdot g_2(D^-)$. Their sum is
$20$ by Proposition \Hgga. This gives that both sides are 10. The proof of the
rest is similar.\qed
\smallskip
For the proofs we use computer calculations. For this one needs a base field.
We take $K=\qz(\zeta)$ where $\zeta$ is an 8th root of unity. The varieties $D^\pm$
are irreducible over $K$. We think that they are irreducible over $\cz$ but it is
not necessary to know this.
\proclaim
{Proposition}
{Let $(g_1,g_2,g_3)$ be three pairwise different elements of $G$. Then
$$g_1(D^\pm)\cdot g_2(D^\pm)\cdot g_3(D^\pm)=2.$$}
Tri%
\finishproclaim
{\it Proof.} 
There have to be treated the following cases. Let $K$ be the subgroup
of $G$ that is generated by $g_1^{-1}g_2$ and $g_1^{-1}g_3$.
\vskip2mm
\item{\rm 1)} The order of $K$ is $16$. There are $16\cdot96$ triples with this property.
There are two subcases:
\itemitem{\rm a)} In $16\cdot48$ cases the intersection of 
$g_1(D^\pm)$, $g_2(D^\pm)$,  $g_3(D^\pm)$ consists of two points 
which are not nodes and which have multiplicity one and which both are defined
over $K$.
\itemitem{\rm b)} In $16\cdot48$ cases the intersection of 
$g_1(D^\pm)$, $g_2(D^\pm)$,  $g_3(D^\pm)$ consists of two points 
which are not nodes and which have multiplicity one and which both are {\it not\/} defined
over $K$.
\vskip1mm\item{\rm 2)} The order of  $K$ is $8$.
There are two subcases:
\itemitem{\rm a)} In $16\cdot36$ cases there are two intersection points
which are defined over $\qz(\imag)$.
None of them is a node and the multiplicity is in each case $1$.
\itemitem{\rm b)} In $16\cdot36$ cases there are $4$ intersection points.
They are all nodes. Here one has to take the intersection numbers in 
the small resolution $\tilde\calX$. In two of the four one has no intersection
point. In each of the other two there is one intersection point with
multiplicity one. This case will be treated in detail in the following Section.
\vskip1mm\item{\rm 3)} The order of $K$ is $4$.
There are $16\cdot42$ cases. In all of them there are $2$ intersection points
which are defined over $\qz(\imag)$. They are not nodes and they have multiplicity $1$.
\smallni
In each of the cases only one sign combination, for example three plus signs,
have to be treated. The other combinations then one can obtain with the help of
Proposition \Hgg.
The only really difficult case is 2b). In all other cases we have intersection
points only outside of the nodes and we can work with the variety $\calX$ and have not
to desingularize it. As in Sect.~2, the intersection points can be 
computed in each case explicitly
and after that one can compute their multiplicities. The case 2b) is more involved.
We will study this case in the next section.
\neupara{Intersections in exceptional fibres}%
In this section we explain the case 2b) in the proof of Proposition \Tri, 
the only case where nodes come into the game.
\smallskip
Before we start, we have to recall some details about the resolutions of a node. 
First, one can consider the blow up of a node. The exceptional fibre is 
biholomorphic equivalent to a
$P^1\times P^1$. 
After a biholomorphic map has been  chosen, we can talk about
horizontal lines $P^1\times \{a_2\}$ and vertical lines $\{a_1\}\times P^1$.
A different choice of the biholomorphic map can preserve ``horizontal'' and
``vertical'' or exchange them. 
One can contract one of the two types of lines to obtain a small resolution
with exceptional fibre $P^1$. Hence there exist two essentially different small
resolutions. The choice of the small resolution
is called a ruling of the node. 
We recall the following fact. Let $Y\subset\calX$ be a pure 
two dimensional subvariety which is smooth at a node. Then the blow-up
of $Y$ gives one of the two small resolutions. The strict transform in the
blow up of the node is a horizontal or vertical line.
\smallskip
The following Lemma has been communicated to us by S.~Cynk.
\proclaim
{Lemma}
{Assume that $a$ is a node in $\calX$ and that $Y_1,Y_2$ are two smooth surfaces which contain $a$.
In the case that the node $a$ is an (maybe embedded)  component of the scheme theoretic
intersection $Y_1\cap Y_2$,
the corresponding lines in $P^1\times P^1$ must be parallel (including equal), i.e.~both
horizontal, or both vertical. 
In the other cases the two lines intersect properly (one horizontal, one vertical).
\medni
{\bf Additional Remark.} In the case of parallel lines, the two lines are equal if and only
if the tangent planes of $Y_1,Y_2$ at $a$ are the same.
}
TwoS%
\finishproclaim
Assume that we have three surfaces $Y_1,Y_2,Y_3$ which have the node $a$ as isolated intersection
point and which are smooth at $a$. 
We consider their strict transforms in the exceptional fibre $P^1\times P^1$ of the
blow up of the node. These can be three parallel lines (equal lines are considered to be
parallel) or two parallel lines and a further line intersecting them.
\vskip0pt
\proclaim
{Definition}
{Let $Y_1,Y_2,Y_3$ be three surfaces in $\calX$ which have the node $a$ 
as isolated intersection point and which are smooth at $a$. They are called
in {\emph good position} (at $a$) if their strict transform in the blow up of the node
consists of two different parallel lines and one which intersects them.}
GooP%
\finishproclaim
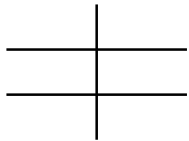
\pspicture[](-4,-1.6)(-2,0.3)
\psline(1,-0.5)(3,-0.5)
\psline(1,-1)(3,-1)
\psline(2,0)(2,-1.5)
\endpspicture
\smallni
We contract this figure in horizontal and vertical direction.
$$\hbox{horizontal}\hbox to 3cm{\hss}\hbox{vertical}
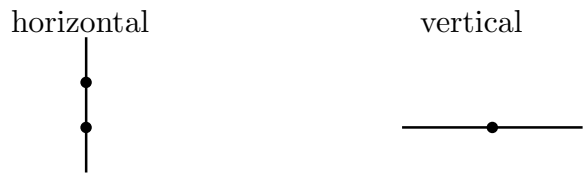
$$
\pspicture[](-4,-0.1)(-2,-1.5)
\psline(0,0)(0,-1.5)
\psline(3.5,-1)(5.5,-1)
\psdots(4.5,-1)(0,-0.5)(0,-1)
\endpspicture
\smallni
In the horizontal direction we see one line and two points on it.
The intersection of these three is empty. In the vertical direction we see
one line and one point on it. The intersection is one point.
Hence we see the following result.
\proclaim
{Proposition}
{Assume that the node $a$ is ruled (so horizontal and vertical is defined).
Assume that we have two triples of surfaces $(D_1',D_2'.D_3')$, 
$(D_1'',D_2'',D_3'')$ which both are in good position in the sense of Definition \GooP.
Assume that their two diagrams are different (not parallel) in the sense of the following figure.
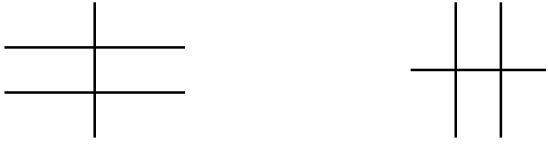
\smallni
\pspicture[](-1,-1.6)(1,0.3)
\psline(1,-0.5)(3,-0.5)
\psline(1,-1)(3,-1)
\psline(2,0)(2,-1.5)
\psline(6,0)(6,-1.5)
\psline(6.5,0)(6.5,-1.5)
\psline(5.5,-0.75)(7,-0.75)
\endpspicture
\smallni
Let $\alpha'$ be the number of intersection points of $D_1',D_2',D_3'$ in the exceptional
$P^1$ and similarly
$\alpha''$. Then
$$\alpha'+\alpha''=1.$$
In particular, this sum is independent of the choice of the ruling.}
IofR%
\finishproclaim 
We consider an example for Proposition \Tri, 2b):
$$\eqalign{
g_1:&=\hbox{identity,}\cr
g_2:&=
(-Y_0,Y_1,-Y_2,Y_3,X_0,-X_1,-X_2,X_3),\cr
g_3:&=(\imag Y_3,-\imag Y_2,-Y_1,-Y_0,X_2,-\imag X_3,\imag X_0,X_1).\cr}$$
Recall that $D^+$ is a component of the
hyperplane section $X_0+X_1+X_2+X_3=0$.
The intersection of the three divisors 
$$g_1(D^+),\quad g_2(D^+),\quad g_3(D^+),$$
considered in $\calX$,
consists of 4 nodes. 
$$\eqalign{
a_1:=&[\sqrt2,0,0,-\sqrt2\imag,0,-1,1],\cr
a_2:=&[\sqrt2,0,0,\sqrt2\imag,0,-1,1],\cr
a_3:=&[-\sqrt2,0,0,-\sqrt2,0,-1,1],\cr
a_4:=&[-\sqrt2,0,0,\sqrt2,0,-1,1].\cr
}$$
Unfortunately the three divisors are not in good position
with respect to all nodes $a_1,\dots,a_4$ 
in the sense of Definition \GooP. For example at the
second node $a_2$ the following happens. The curves $g_1(D^+)\cap g_3(D^+)$ and 
$g_2(D^+)\cap g_3(D^+$) have $a_2$ as isolated fixed point. But the intersection
of $g_1(D^+)$ and $g_3(D^+)$ in $\calX$
is set theoretically an elliptic curve,
cut out in $\calX$ by the linear equations
$$Y_1=Y_2=X_0+X_3=X_1+X_2.$$
But there are four embedded components in this curve, namely the $4$ nodes.
Hence the three lines in the exceptional $P^1\times P^1$ of the node are parallel.
\smallskip
To remedy this situation, we replace the divisor $g_1(D^+)$
by the other component $g_1(D^-)$.
From Proposition \Hgg\ we see that
$$g_1(D^+)\cdot g_2(D^+)\cdot g_3(D^+)+g_1(D^-)\cdot g_2(D^+)\cdot g_3(D^+)=4.$$
Hence we see
$$g_1(D^+)\cdot g_2(D^+)\cdot g_3(D^+)=2\Longleftrightarrow 
g_1(D^-)\cdot g_2(D^+)\cdot g_3(D^+)=2.$$
\proclaim
{Remark} 
{We use now the notation
$$D_1=g_1(D^-),\quad D_2=g_2(D^+),\quad D_3=g_3(D^+).$$
\item{\rm a)} The intersection of $D_1\cap D_2\cap D_3$ consists of the same $4$
nodes $a_1,\dots,a_4$.
\item{\rm b)} The three divisors $D_1,D_2,D_3$ are in good position at all $4$
nodes.}
Dsas%
\finishproclaim
We can use Lemma \TwoS\ to clarify the position of two of the $D_i$. The computation can be done
with help of MAGMA which allows to compute the primary decomposition of an ideal.\qed
\smallskip
We need now some information about the small desingularization $\tilde\calX$.
The group $G$ has the essential fact that it extends to $\tilde\calX$. This means that it
preserves the rulings of the nodes. We use this to compare the rulings of the nodes
$a_1,a_2,a_3,a_4$. One checks
$$g_1(a_2)=a_1,\quad g_1(a_4)=a_3.$$
We want to compare the number of intersection points of $D_1,D_2,D_3$ over the node $a_1$ with 
those over $a_2$. Instead of this we can compare 
\vskip1pt
\item{\rm a)} the number of intersection points of $D_1,D_2,D_3$ over the node $a_1$,
\item{\rm b)} the number of intersection points of $g_1(D_1),g_1(D_2),g_1(D_2)$ over $a_1$.
\smallni
One can compute the positions of the line diagrams in the sense of Proposition \IofR.
It turns out that the two diagrams are different. So we can apply Proposition \IofR\ to
obtain the following result.
\smallskip
The divisors $D_1,D_2,D_3$ have an intersection point in the small resolution over $a_1$ or
$a_2$ but not over both. The same argument works for the pair $a_3,a_4$. Hence we get
$$D_1\cdot D_2\cdot D_3=2.$$
We settled an example for 2b). Each other case of Proposition \Tri.
can be treated in the same manner.
This completes the proof of Proposition \Tri, 2b).\qed
\neupara{Intersection numbers in the quotient}%
We recall from the paper [FS].
\proclaim
{Theorem}
{The orbit of the subvariety $D^+$ (the same is true for $D^-$) under the group $G$ consists of
$16$ subvarieties of $\calX$. They are non-singular.  Each of the 96
nodes is contained in 4 of the 16 subvarieties.
The blow up along the union of the 16 (considered as reduced subvariety)
is a smooth projective variety $\tilde\calX$ with a free $G$-action. The quotient
$\tilde\calX/G$ is a (projective) rigid Calabi-Yau manifold with $h^{11}=2$ ($e=4$).
}
diV%
\finishproclaim
From [FS] we know the following result.
\proclaim
{Proposition}
{We denote by $\calD^\pm$ the image of $D^\pm$ in $\tilde\calX/G$.
The two divisors give a $\qz$-basis of
$$\Pic(\tilde\calX/G)\otimes_\gz\qz.$$}
baS%
\finishproclaim
We want to compute the intersection numbers $\calD^\pm\cdot\calD^\pm\cdot\calD^\pm$.
These are essentially $4$ cases, depending on the number of the plus  (or minus) signs.
Let $\pi:\calX\to\calX/G$ be the canonical projection. We use the formula
$$16\calD^\pm\cdot\calD^\pm\cdot\calD^\pm=
\pi_*(\pi^*\calD^\pm\cdot\pi^*\calD^\pm\cdot\pi^*\calD^\pm).$$
It implies
$$16\calD^\pm\cdot\calD^\pm\cdot\calD^\pm=
\pi^*\calD^\pm\cdot\pi^*\calD^\pm\cdot\pi^*\calD^\pm.$$
We have
$$\pi^*\calD^\pm=\sum_{g\in G} g(D^\pm).$$
We get
$$\eqalign{
\calD^\pm\cdot\calD^\pm\cdot\calD^\pm=&{1\over 16}
\sum_{g_1,g_2,g_3}g_1(D^\pm)\cdot g_2(D^\pm)\cdot g_3(D^\pm)=\cr&
\sum_{g_1,g_2}D^\pm\cdot g_1(D^\pm)\cdot g_2(D^\pm).\cr}$$
With the results of the previous section we can derive now our main result.
\proclaim
{Theorem}
{The two divisors $\calD^\pm$ of the rigid Calabi-Yau manifold $\tilde\calY:=\tilde\calX/G$
give a basis of $\Pic\tilde\calY\otimes_\gz\qz$. Their intersection numbers are
$$\calD^+\cdot\calD^+\cdot\calD^+=\calD^-\cdot\calD^-\cdot\calD^-=2^7,\quad
\calD^+\cdot\calD^+\cdot\calD^-=\calD^+\cdot\calD^-\cdot\calD^-=5\cdot 2^7.$$
}
InterS%
\finishproclaim
The intersection numbers get nicer if one uses a different basis of $\Pic\calY\otimes_\gz\qz$,
namely
$$\calA={1\over 16}(\calD^++\calD^-),\quad \calB={1\over 8}(\calD^+-\calD^-).$$
Then we get
\rahmen{\calA^3=1,\quad \calB^3=\calA^2\calB=0,\quad \calA\calB^2=-1.}
There may be several projective small resolutions $\tilde\calX$ of $\calX$
such that $G$ extends to it. The results show that the intersection numbers are
in all cases the same.
\vskip1.5cm\noindent
{\paragratit References}
\bigskip
\item{[CM]} Cynk, S., Meyer, C.:
{\it Modular Calabi-Yau Threefolds of level eight,\/}
Internat.~J.~Math.\ {\bf 18}, no.~3, 331--347 (2007)
\medni
\item{[FS]} Freitag, E., Salvati--Manni, R.:
{\it On Siegel threefolds with a projective Calabi--Yau model,\/}
Commun. Number Theory Phys.\ 5, No.\ 3, 713--750 (2011) 
\medskip
\item{[GN]} van Geemen, B., Nygaard, N.O.:
{\it On the geometry and arithmetic of some Siegel modular threefolds,\/}
Journal of Number Theory {\bf 53}, 45--87 (1995)
\medskip
\item{[Sh]} Shafarevich, I.R.: {\it Basic Algebraic Geometry,}
Springer--Verlag Berin Heidelberg New York (1974)
\medni
{\it Eberhard Freitag\hfill\break
Mathematisches Institut\hfill\break
Im Neuenheimer Feld 288\hfill\break
D69120 Heidelberg}\hfill\break
{\tt freitag@mathi.uni-heidelberg.de}
\bye